\def\cO{{\cal O}}

\documentclass{amsart}

\usepackage{amsmath, amscd, amssymb}
\usepackage[all]{xy}
\usepackage{color}

\newcommand{\cal}{\mathcal}

\newcommand{\bC}{{\Bbb C}}

\newcommand{\bF}{{\Bbb F}}

\newcommand{\bP}{{\Bbb P}}

\newcommand{\bV}{{\Bbb V}}
\newcommand{\bZ}{{\Bbb Z}}

\newcommand{\cP}{{\cal P}}

\newcommand{\cW}{{\cal W}}

 \DeclareMathOperator{\ch}{ch}
\DeclareMathOperator{\reg}{reg}

\DeclareMathOperator{\rank}{rank} 
\DeclareMathOperator{\td}{td} \DeclareMathOperator{\tr}{tr}

\newtheorem{theorem}{Theorem}[section]

\newtheorem{proposition}[theorem]{Proposition}
\newtheorem{lemma}[theorem]{Lemma}

\theoremstyle{remark}
\newtheorem{remark}{Remark}[section]

\theoremstyle{definition}

\newtheorem{definition}{Definition}[section]

\definecolor{light-pink}{rgb}{1,.90,.90}
\definecolor{light-green}{rgb}{.95,1,.95}

\definecolor{yellow}{rgb}{1,1,0}
\definecolor{orange}{rgb}{1,.7,0}
\definecolor{red}{rgb}{1,0,0}
\definecolor{white}{rgb}{1,1,1}

\definecolor{A}{rgb}{.75,1,.75}

\begin{document}

\title{Topological String Partition Functions as Equivariant Indices}

\author{Jun Li}
\address{Department of Mathematics\\
Stanford University, Stanford, CA, USA}
\email{jli@math.stanford.edu}

\author{Kefeng Liu}
\address{Center of Mathematical Sciences, Zhejiang University, Hangzhou,
China and Department of Mathematics\\
University of California at Los Angeles} \email{liu@math.ucla.edu,
liu@cms.zju.edu.cn}

\author{Jian Zhou}
\address{Department of Mathematical Sciences\\
Tsinghua University\\ Beijing, China}
\email{jzhou@math.tsinghua.edu.cn}

\maketitle

\begin{abstract}
We propose to use the identification of topological string
partition functions as equivariant indices on framed moduli spaces
of instantons to study the Gopakumar-Vafa conjecture for some
local Calabi-Yau geometries.
\end{abstract}


\section{Introduction}

In this work we exploit the relationship with certain equivariant genera of isntanton moduli spaces
to study the string partition functions of some
local Calabi-Yau geometries,
in particular,
the Gopakumar-Vafa conjecture for them \cite{Gop-Vaf}.

Gromov-Witten invariants are in general rational numbers. However
as conjectured by Gopakumar and Vafa \cite{Gop-Vaf} using
M-theory, the generating series of Gromov-Witten invariants in all
degrees and all genera has a particular form, determined by some
integers. See \cite{Gop-Vaf} or Section \ref{sec:Prelim} for
precise formulation. There have been various proposals to the
proof of this conjecture \cite{Kat-Kle-Vaf, Hos-Sai-Tak, MNOP}.
 In this work we propose a new and geometric
approach towards this conjecture for some interesting cases. The
method relates the computations of Gromov-Witten invariants to
equivariant index theory and 4 dimensional gauge theory.

Recently there have been some progresses on the calculations of
Gromov-Witten invariants, both in the physical approaches
\cite{Aga-Mar-Vaf, Iqb, Aga-Kle-Mar-Vaf} and mathematical
treatments \cite{ZhoFree, LLLZ}. For toric local Calabi-Yau
geometries, one can now express the string partition functions as
sums over partitions. Furthermore, one can relate them to
partition functions in gauge theory according to the idea of
geometric engineering (see e.g. \cite{Nek, Iqb-Kas1, Iqb-Kas2,
Egu-Kan1, ZhoCounting, Egu-Kan2, Hol-Iqb-Vaf}). The latter are
equivariant genera of framed moduli spaces and can be computed by
localization formula \cite{Nak-Yos, Hol-Iqb-Vaf} (see also Section
\ref{sec:Framed}). The fixed points on moduli spaces are tuples of
partitions hence one gets seemingly different sums over
partitions. However as shown in some of the works mentioned above,
one can use combinatorics to identify different expressions.

To prove the GV conjecture,
one needs to rewrite the sums over partitions
as infinite products (see e.g. \cite{Hol-Iqb-Vaf} or Section \ref{sec:Prelim} for explanation).
The purpose of this paper is to point out that in some cases if one pushes forward
the calculations from the framed moduli spaces,
which are the Gieseker partial compactification of the moduli space of genuine instantons,
to the Uhlenbeck partial compactification,
then one can achieve this.
Let us briefly explain some relevant terminologies.

For a compact complex $d$-manifold $X$,
one is often interested in its Hirzebruch $\chi_y$ genus.
It is defined by:
$$\chi_y(X) = \sum_{p=0}^d (-y)^p \sum_{q=0}^d (-1)^q \dim H^q(X, \Lambda^pT^*X).$$
By Hirzebruch-Riemann-Roch theorem, it can be computed as follows:
\begin{eqnarray*}
&& \chi_y(X) =  \int_X \prod_{j=1}^d \frac{x_j(1- ye^{-x_j})}{(1- e^{-x_j})},
\end{eqnarray*}
where $x_1, \dots, x_d$ are the formal Chern roots of $TX$.
The $\chi_y$-genus reduces to other invariants for special values of $y$,
e.g., $\chi_1(X)$ is the Euler number,
and
$$\chi_0(X) = \sum_{q=0}^d (-1)^q \dim H^q(X, \cO_X)$$
is the geometric genus of $X$.
An important generalization of the $\chi_y$-genus is the elliptic genus $\chi(X, y, q)$,
which can be defined as the generating series of dimensions of some cohomology
groups of series of vector bundles.
It can be computed by
\begin{eqnarray*}
&& \chi(X, y, q) = y^{-d/2} \int_X \prod_{j=1}^d x_j
\prod_{n \geq 1}\frac{(1-yq^{n-1}e^{-x_j})(1- y^{-1}q^ne^{x_j})}
{(1-q^{n-1}e^{-x_j})(1-q^ne^{x_j})}.
\end{eqnarray*}
It is easy to see that
$$\chi_y(X) = y^{d/2} \chi(X, y, 0).$$

Denote by $M(N, k)$ the framed moduli space of torsionfree sheaves
of rank $N$ and $c_2 = k$ on  $\bP^2$. The instanton partition
functions in 4D, 5D and 6D correspond to the generating series of
the $\chi_0$, $\chi_y$, and the elliptic genera, respectively, of
$M(N, k)$. Note however $M(N, k)$ is noncompact hence the
cohomology groups might be infinite dimensional, so the definition
of the genera need to be modified. It turns out that the weight
spaces of the cohomology groups with respect to some natural torus
actions are finite-dimensional, hence it makes sense to replace
the dimensions of the cohomology groups by their characters. In
this way, one gets the equivariant versions of the $\chi_0$,
$\chi_y$ and elliptic genera of $M(N, k)$. Furthermore, one can
compute them by using localization formula.

Since the fixed points on $M(N, k)$ are parameterized by tuples of partitions,
localization methods applied in this gauge theory setting yield sums over tuples of partitions.
On the string theory side,
one has different looking expressions as sums over partitions.
However combinatorics can be used to identity these expressions
\cite{Iqb-Kas1, Iqb-Kas2, Egu-Kan1, ZhoCounting, Egu-Kan2}.
Hence one can relate equivariant genera of framed moduli spaces in gauge theory
to the generating series of the Gromov-Witten
invariants of some local toric Calabi-Yau geometries.
This idea appeared in physics literature in the context of geometric engineering,
for example Nekrasov's conjecture \cite{Nek}.
Our detailed treatment includes the cases with matters,
hence it makes the geometric engineering of $SU(N)$ gauge theory with matters
mathematically rigorous.

We notice that the identification of the Gromov-Witten partition
functions as equivariant genera suggests an approach to prove the
Gopakumar-Vafa conjecture for the relevant Calabi-Yau geometries.
To prove this conjecture, one needs to rewrite the sums over
partition as infinite products. We propose a natural geometric
setting to achieve this. More precisely, by exploiting the natural
morphism $M(N, k) \to M(N, k)_0$ from the framed moduli spaces
$M(N, k)$ to the Uhlenbeck compactifications $M(N, k)_0$, one can
pushforward the calculations of equivariant indices on $M(N, k)$
to $M(N, k)_0$. When $N=1$, $M(1, k)$ is the Hilbert scheme
$(\bC^2)^{[k]}$ while $M(1, k)_0$ is the symmetric product
$(\bC^2)^{(k)}$. By standard manipulation on the symmetric
products (cf. e.g. \cite{Zho991, Zho992}) one can get infinite
product expressions and prove the Gopakumar-Vafa conjecture. We
will discuss in detail the case of $\chi_0$, $\chi_y$, and
elliptic genera which correspond to three different local
Calabi-Yau geometries. We propose to deal with the $N \geq 2$
cases in a similar fashion.

We summarize the steps involved in our approach as follows.
\begin{itemize}
\item[Step 1] Compute the Gromov-Witten invariants by diagrammatic methods as sum over partitions;

\item[Step 2] Use Schur function calculus to rewrite the result in Step 1 as
sum over partitions of certain ratios;

\item[Step 3] Identify the result in Step 2 as suitable equivariant genera on $M(N, k)$;

\item[Step 4] Identify the result in Step 3 as suitable equivariant genera on $M_0(N, k)$;

\item[Step 5] Obtain Gopakumar-Vafa type infinite products by series manipulations.
\end{itemize}
The first three steps are mostly from existing works,
and the last two steps are the new features of this work.

Note it is a common practice in number theory that if one wants to show certain number is integral,
then one tires to identify it as the dimension of certain cohomology group.
Also it is well-known that by taking symmetric powers,
one can identify a sum expression with a product expression.
Our proposal is compatible with such general principles.

The rest of the paper is arranged as follows. In Section
\ref{sec:Prelim} we recall the Gopakumar-Vafa Conjecture and its
infinite product formulation. In Section \ref{sec:Symm} we compute
the equivariant elliptic genera for symmetric products, which
naturally yield some infinite product expressions. We recall some
results on framed moduli spaces and compute their equivariant
ellitpic genera in Section \ref{sec:Framed}. The relationship
between the equivariant genera of Hilbert schemes and symmetric
products is studied in Section \ref{sec:Hilbert}, and applied in
Section \ref{sec:N=1} to compute some Gopakumar-Vafa invariants
for Calabi-Yau geometries corresponding to $M(1, k)$. We give two
theorems about the index expressions of the Gromov-Witten
invariants of certain toric Calabi-Yau manifolds, and propose a
similar approach of infinite products for Calabi-Yau geometries
corresponding to $M(N, k)$ for $N > 1$ by push-froward in
equivariant $K$-theory. The detail will be worked out later.

The major part of this paper was done during the authors' visit of
Center of Mathematical Science, Zhejiang University in 2004. The
authors would like to thank CMS for its hospitality. The first and
the second author are supported by NSF and the thrid author by
NSFC. During the preparation of the paper, there have been
interesting progresses in the subject. See \cite{Pan} and
\cite{Kon1} for the proof of the GV conjecture for toric
Calabi-Yau manifolds by using combinatorics and the theory of
topological vertex.


\section{Preliminaries}
\label{sec:Prelim}

In this section we recall the Gopakumar-Vafa Conjecture and its formulation
in terms of infinite product.

\subsection{Gopakumar-Vafa Conjecture}

Denote by $F_X$ the generating series of Gromov-Witten invariants of a Calabi-Yau
$3$-fold $X$.
Intuitively,
one counts the number of stable maps with {\em connected} domain curves to $X$
in any given nonzero homology classes.
However, because the existence of automorphisms,
one has to performed the weighted count by dividing by the order of automorphism groups.
Hence Gromov-Witten invariants are in general rational numbers.

Based on $M$-theory considerations, Gopakumar and Vafa \cite{Gop-Vaf} made a remarkable
conjecture on the structure of $F$,
in particular, its integral properties.
More precisely,
there are integers $n^g_{\Sigma}$ such that
\begin{eqnarray}
F & = & \sum_{\Sigma \in H_2(X)-\{0\}} \sum_{g \geq 0} \sum_{k \geq 1} \frac{1}{k}
n^g_{\Sigma} (2\sin \frac{k\lambda}{2})^{2g-2}
Q^{k\Sigma}.
\end{eqnarray}
Let $q = e^{\sqrt{-1}\lambda}$,
and regard it as an element of $SU(2)$ represented by the diagonal matrix
$$\begin{pmatrix} q & 0 \\ 0 & q^{-1} \end{pmatrix}.$$
Denote by $V_n$ the $(n+1)$-dimensional irreducible representation of $SU(2)$.
Then the character of $V_n$ is given by:
$$\chi_{V_n}(q) = \frac{q^{n+1} - q^{-(n+1)}}{q - q^{-1}}
= q^{n} + q^{n-2} + \cdots + q^{-n}.$$
We have
$$(2\sin \frac{\lambda}{2})^2 = - (q+q^{-1}+2)
= -(\chi_{V_1} + 2\chi_{V_0}).$$
Recall $\{V_n: n \geq 0\}$ form an integral basis
of the representation ring of $SU(2)$.
Since $V_1 \otimes V_n = V_{n+1} \oplus V_{n-1}$,
$\{(V_1 \oplus V_0 \oplus V_0)^{\otimes n}: n \geq 0\}$ also form an integral basis.
Therefore,
there are integers $N_{\Sigma}^g$ such that
\begin{eqnarray} \label{eqn:nN}
\sum_{g \geq 0} n^g_{\Sigma} (-1)^g (q^{\frac{1}{2}} - q^{-\frac{1}{2}})^{2g}
= \sum_{g \geq 0}N_{\Sigma}^g (q^{g} + q^{g-2} + \cdots + q^{-g}).
\end{eqnarray}

\subsection{Infinite product formulation of the Gopakumar-Vafa Conjecture}
The generating series of {\em disconnected} Gromov-Witten invariants are given by
the string partition function:
$$Z = \exp F.$$

See \cite{Hol-Iqb-Vaf} for the following:

\begin{proposition}
The Gopakumar-Vafa Conjecture can be reformulated as follows:
\begin{eqnarray} \label{eqn:GV}
&& Z = \prod_{\Sigma \in H_2(X)}\prod_j \prod_{k=-j}^j \prod_{m=0}^{\infty}
(1- q^{2k+m+1}Q^{\Sigma})^{(-1)^{2j+1}(m+1)N_{\Sigma}^g}.
\end{eqnarray}
where $j = \frac{g}{2}$,
$k=-j, -j+1, \dots, j-1, j$.
\end{proposition}

\begin{proof}
By (\ref{eqn:nN}) one has
\begin{eqnarray*}
&& \sum_{\Sigma \in H_2(X)} \sum_{g \geq 0} \sum_{k \geq 1} \frac{1}{k}
n^g_{\Sigma} (2\sin \frac{k\lambda}{2})^{2g-2}
Q^{k\Sigma} \\
& = & \sum_{\Sigma \in H_2(X)} \sum_{g \geq 0} \sum_{k \geq 1} \frac{(-1)^{g-1}}{k}
N^g_{\Sigma} \frac{\sum_{a=0}^{g} q^{k(g-2a)}}{(q^{\frac{k}{2}} - q^{-\frac{k}{2}})^2}
Q^{k\Sigma} \\
& = & \sum_{\Sigma \in H_2(X)} \sum_{g \geq 0} \sum_{k \geq 1} \frac{(-1)^g{-1}}{k}
N^g_{\Sigma} \frac{\sum_{a=0}^{g} q^{k(g-2a+1)}}{(1 - q^k)^2}
Q^{k\Sigma} \\
& = & \sum_{\Sigma \in H_2(X)} \sum_{g \geq 0} \sum_{k \geq 1} \frac{(-1)^{g-1}}{k}
N^g_{\Sigma} \sum_{a=0}^{g} q^{k(g-2a+1)}
Q^{k\Sigma} \sum_{m \geq 0} (m+1) q^{km} \\
& = & \sum_{\Sigma \in H_2(X)} \sum_{g \geq 0} \sum_{k \geq 1} \frac{(-1)^{g-1}}{k}
N^g_{\Sigma} \sum_{a=0}^{g} q^{k(g-2a+1)}
Q^{k\Sigma} \sum_{m \geq 0} (m+1) q^{km} \\
& = & -\sum_{\Sigma \in H_2(X)} \sum_{g \geq 0} (-1)^g N^g_{\Sigma} \sum_{a=0}^{g} \sum_{m \geq 0} (m+1)
\sum_{k \geq 1} \frac{1}{k} (q^{g-2a+m+1} Q^{\Sigma})^k \\
& = & \sum_{\Sigma \in H_2(X)} \sum_{g \geq 0} (-1)^g N^g_{\Sigma} \sum_{a=0}^{g} \sum_{m \geq 0} (m+1)
\log (1- q^{g-2a+m+1} Q^{\Sigma}).
\end{eqnarray*}
\end{proof}

>From (\ref{eqn:GV}) one sees that to prove the Gopakumar-Vafa conjecture one needs
an infinite product expression for the string partition function.
On the other hand,
as will be recalled below
such partition functions are usually given by a sum expression in the diagrammatic method.
We will show by some examples how one can convert the sum expressions
to the product expressions.
This is achieved by relating the string partition functions
to equivariant genera, first of Hilbert schemes, then of symmetric products.

Our motivation is very simple.
As noted in \cite{Hol-Iqb-Vaf},
the expression on the right-hand side of (\ref{eqn:GV})
looks like `counting' the states in the Hilbert space of
a second quantized theory,
a situation similar to the calculations of
orbifold elliptic genera of symmetric products \cite{Dij-Moo-Ver-Ver}.


\section{Orbifold Equivariant Elliptic Genera of Symmetric Products}
\label{sec:Symm}

In this section we show that infinite product expressions arise naturally when
one considers genera of symmetric products.
This motivates our proposal of using symmetric products to prove the Gopakumar-Vafa Conjecture.
For references see e.g. \cite{Zho991, Zho992, Bor-Lib}.

\subsection{Orbifold elliptic genera}

Let $M$ be a compact complex manifold and let $G$ act on $M$ by biholomorphic maps.
Furthermore,
let $\pi: E \to M$ be a holomorphic vector bundle which admits a $G$-action compatible with the
$G$-action on $M$.
For $g \in G$,
the Lefschetz number is defined by:
$$L(M, E)(g) = \sum_{i=0}^{\dim M} (-1)^i \tr (g|_{H^i(M, \cO(E))}).$$
It can be computed by the holomorphic Lefschetz formula \cite{Ati-Sin}:
$$L(M, E)(g) = \int_{M^g} \frac{\ch_g E}{\ch_g \Lambda_{-1} (N_{M^g/M})},$$
where $M^g$ denotes the set of points fixed by $g$,
$N_{M^g/M}$ denotes the normal bundle of $M^g$ in $M$.
Here and in the following it is understood that one considers each connected component
of $M^g$ separately.

There is a natural decomposition
$$TM|_{M^g} = \oplus_j N^g(\lambda_j),$$
where each $N^g(\lambda_j)$ is a holomorphic subbundle on which $g$ acts as
$e^{2\pi \sqrt{-1} \lambda_j}$
where $0 \leq \lambda_j < 1$.
In particular,
$$N^g(0) = TM^g.$$
Define the {\em fermionic shift} by \cite{Zas}:
\begin{eqnarray} \label{shift}
F(M^g) = \sum_j (\rank N_{\lambda_j}) \lambda_j.
\end{eqnarray}
Following \cite{Bor-Lib},
define
\begin{eqnarray*}
E(M;q,y)^g & = & y^{- \frac{d}{2} + F(M^g_i)}
\otimes_j \left[ \otimes_{n \geq 1} \left(
\Lambda_{-yq^{n - 1}}(TM^g)^*
\otimes \Lambda_{-y^{-1} q^{n}} TM^g
\right) \right. \\
&& \left.  \otimes \otimes_{n \geq 1}
\left( S_{q^n} (TM^g)^* \otimes S_{q^n} TM^g \right) \right] \\
&& \otimes_{\lambda_j \neq 0} \left[ \otimes_{n \geq 1} \left(
\Lambda_{-yq^{n - 1 + \lambda_j}}N_{\lambda_j}^*
\otimes \Lambda_{-y^{-1} q^{n - \lambda_j}} N_{\lambda_j}
\right) \right. \\
&& \left.  \otimes \otimes_{n \geq 1}
\left( S_{q^{n -1 + \lambda_j}} N_{\lambda_j}^*
\otimes S_{q^{n - \lambda_j}} N_{\lambda_j} \right) \right].
\end{eqnarray*}
Recall where for complex vector bundle $E$,
\begin{align*}
& \Lambda_t(E) = 1 + tE + t^2\Lambda^2E + \cdots, \\
& S_t(E) =1+tE+t^2 S^2E+\cdots .
\end{align*}

\begin{definition}
The orbifold elliptic genus of $M/G$ is defined by
$$\chi(M, G; q, y)
=\sum_{[g] \in G_*} \frac{1}{|Z(g)|}
\sum_{h \in Z(g)} L(M^g, E(M; q, y)^g)(h),$$
where $G_*$ denotes the set of conjugacy classes of $G$,
and $Z(g)$ denotes the centralizer of $g$.
\end{definition}

For any $h \in Z(g)$,
let $M^{(g, h)}$ be the set of points of $M$
fixed by both $g$ and $h$.
By Lefschetz formula,
we have
\begin{eqnarray*}
&& L(M^g, E(M; q, y)^g)(h)
= \int_{M^{g, h}} \frac{\ch_h(E(M; q, y)^g)}
{\ch_h(\Lambda_{-1}(N^*_{M^{g,h}/M^g}))}
\td(TM^{g, h}).
\end{eqnarray*}

\subsection{Equivariant orbifold elliptic genera}

Assume now a $G$-manifold $M$
admits an $S^1$-action which commutes with the $G$-action.
The equivariant orbifold elliptic genus is defined by:
$$\chi(M, G; q, y)(s) = \sum_{[g] \in G_*}
\frac{1}{|Z(g)|} \sum_{h \in Z(g)}
L(M^g, E(M; q, y)^g)(h, s),$$
where $s \in S^1$.
This defines a character for $S^1$.
Denote by $M^{(g, h, s)}$ the points on $M$ which is fixed by
$g$, $h$, and $s$.
For our applications,
it suffices to assume that $M^{g, h, s}$ consists of an isolated point.
Then by Lefschetz formula,
we have
\begin{eqnarray*}
&& L(M^g, E(M; q, y)^g)(h, s)
=   \frac{\ch_{h, s}(E(M; q, y)^g)}
{\ch_{h, s}(\Lambda_{-1}(N^*_{M^{g,h,s}/M^g}))}.
\end{eqnarray*}

\subsection{Equivariant elliptic genera of symmetric products}

Suppose $X$ is a nonsingular projective variety admitting an action by a torus group $T$.
The diagonal action by $T$ and the natural action by permutation group $S_N$
on the $N$-fold cartesian product $X^N$ commute with each other.

\begin{theorem} \label{thm:Symmetric}
Suppose the equivariant elliptic genera of $X$ can be written as
$$\chi(X; q, y)(t_1, \dots, t_r)
= \sum_{m \geq 0, l, k} c(m, l, k) q^my^lt_1^{k_1} \cdots t_r^{k_r},$$
then one has
\begin{equation} \label{eqn:SymmEll}
\begin{split}
&\sum_{N \geq 0}  Q^N\chi(X^N, S_N; q, y)(t_1, \dots, t_r) \\
= & \exp \sum_{N, n \geq 1} \frac{Q^{Nn}}{N}
\cdot \frac{1}{n} \sum_{i=0}^{n-1} \chi(X)((\omega_n^iq^{1/n})^{N}, y^{N})(t_1^{N}, \cdots, t_r^{N}) \\
= & \prod_{n > 0, m \geq 0, l, k_1, \dots, t_r}
\frac{1}{(1-Q^nq^my^lt_1^{k_1} \cdots t_r^{k_r})^{c(nm, l, k_1, \dots, k_r)}},
\end{split} \end{equation}
where $\omega_n = \exp (2\pi\sqrt{-1}/n)$.
\end{theorem}

Note we have
\begin{eqnarray*}
&& \sum_{N \geq 0}  Q^N\chi_y(X^N, S_N)(t_1, \dots, t_r)
= \sum_{N \geq 0} (yQ)^N\chi(X^N, S_N; 0, y) (t_1, \dots, t_r).
\end{eqnarray*}
Hence by taking $q=0$ in (\ref{eqn:SymmEll}) one gets:
\begin{equation} \label{eqn:SymmChiy}
\begin{split}
& \sum_{N \geq 0}  Q^N\chi_y(X^N, S_N)(t_1, \dots, t_r) \\
= & \exp \sum_{m \geq 1} \frac{Q^m\chi_{y^m}(X)(t_1^m, \dots, t_r^m)}{m(1-y^mQ^m)} \\
= & \prod_{n > 0,  l, k_1, \dots, k_r}
\frac{1}{(1-(yQ)^ny^lt^k)^{c(0, l, k_1, \dots, k_r)}}.
\end{split} \end{equation}
Taking $y=0$ in (\ref{eqn:SymmChiy}) one gets:
\begin{eqnarray} \label{eqn:SymmChi0}
&& \sum_{N \geq 0}  Q^N\chi_0(X^N, S_N)(t_1, \dots, t_r)
= \exp (\sum_{m \geq 1} \frac{Q^m}{m}  \chi_{0}(X)(t_1^m, \dots, t_r^m)).
\end{eqnarray}
See \cite{Zho991, Zho992, Bor-Lib} for the nonequivariant version of
(\ref{eqn:SymmEll})-(\ref{eqn:SymmChi0}).

\subsection{Equivariant elliptic genera of the symmetric products: The proof}
Assume
$$TX|_{X^s} = \oplus_v N_v,$$
where $s = e^{2\pi \sqrt{-1} t}$ acts on $N_v$
by multiplication by $e^{2\pi \sqrt{-1} v t}$.
Let $g = m\sigma_l$,
$h = h_{l,j, m}$, and $n= m l$.
Here $\sigma_l$ denotes an $l$-cycle,
$m\sigma_l$ stands for the composition of $m$ disjoint $l$-cycles.
For examples,
the permutation $(12)(34)(56)$ can be written as $3\sigma_2$.
Then we have
\begin{eqnarray} \label{eqn:ghsweights}
TX^n|_{(X^n)^{g, h,s}} & = & \oplus_{k=0}^{l-1}\oplus_{r=0}^{m - 1}
\oplus_v N_{k, r, v},
\end{eqnarray}
where $N_{k, r, v}$ is a copy of $N_v$ on which
$g$ acts as multiplication by $e^{2\pi \sqrt{-1} k/l}$,
and $h$ acts as multiplication by
$e^{2\pi \sqrt{-1} (r/m + jk/(lm))}$.

\begin{lemma}
If $\chi(X; q, y)(s) = \sum c(\alpha, \beta, \gamma)
q^\alpha y^\beta s^\gamma$,
then
\begin{eqnarray} \label{eqn:Lefschetz3}
&& \frac{1}{l} \sum_{j = 0}^{l-1} L((X^{lm})^{(m \sigma_l)},
    E(X^{lm}; q, y)^{(m \sigma_l)})(h_{l, j, m}, s)
= \sum_{\alpha, \beta, \gamma} c(l \alpha, \beta, \gamma)
q^{m \alpha} y^{m \beta} s^{m \gamma}.
\end{eqnarray}
\end{lemma}

\begin{proof}
Let $d_v = \dim N_v$.
Denote the formal Chern roots of $N_v$ by
$2\pi \sqrt{-1} x_{vi}$,
$i =1, \dots, d_v$.
Then we have
\begin{eqnarray*}
&& \chi(X;q, y)(s) = \sum_{\alpha, \beta, \gamma}
c(\alpha, \beta, \gamma) q^{\alpha}y^{\beta} s^{\gamma} \\
& = & \int_{X^s} \prod_{i=1}^{d_0} (2\pi\sqrt{-1} x_{0i})
\cdot \prod_v \prod_{i=1}^{d_v} \frac{\theta(x_{vi} - z + v t, \tau)}
{\theta(x_i + v t, \tau)}.
\end{eqnarray*}
On the other hand,
\begin{eqnarray*}
&& L((X^{lm})^{(m \sigma_l)},
    E(X^{lm}; q, y)^{(m \sigma_l)})(h_{l, j, m}, s) \\
& = & \int_{X^s} \prod_{i = 1}^{d_0} (2 \pi \sqrt{-1} x_{0i})
\cdot \prod_v \prod_{i=1}^{d_v} \prod_{k=0}^{l-1} \prod_{r=0}^{m-1}
y^{k/l}
\frac{\theta(x_{vi} - z - k\tau/l + r/m + jk/(lm) + v t, \tau)}
{\theta(x_{vi} - k\tau/l + r/m + jk/(lm) + v t, \tau)} \\
& = & \int_{X^s} \prod_{i = 1}^{d_0} (2 \pi \sqrt{-1} x_{0i})
\cdot \prod_v  \prod_{i=1}^{d_v} \prod_{k=0}^{l-1} y^{km/l}
\frac{\theta(m x_{vi} - m z - m k\tau/l + jk/l + m v t, m \tau)}
{\theta(m x_{vi} - m k \tau + jk/l + m v t, m\tau)} \\
& = & \int_{X^s} \prod_{i = 1}^{d_0} (2 \pi \sqrt{-1} x_{0i})
\cdot \prod_v \prod_{i=1}^{d_v}
\frac{\theta(m x_{vi} -m z + m v t, (m \tau - j)/l)}
{\theta(m x_{vi} + m v t, (m\tau - j)/l)} \\
& = & \frac{1}{m^{d_0}} \int_{X^s} \prod_{i = 1}^d
(2 \pi \sqrt{-1} m x_{0i})
\cdot \prod_v  \prod_{i=1}^{d_v}
\frac{\theta(m x_{vi} - m z + m v t, (m \tau - j)/l)}
{\theta(m x_{vi} + m v t, (m\tau - j)/l)} \\
& = & \int_{X^s} \prod_{i = 1}^{d_0} (2 \pi \sqrt{-1} x_{0i})
\cdot \prod_v  \prod_{i=1}^{d_v}
\frac{\theta(x_{vi} - m z + m v t, (m \tau - j)/l)}
{\theta(x_{vi} + m v t, (m\tau - j)/l)} \\
& = & \sum_{\alpha, \beta, \gamma} c(\alpha, \beta, \gamma)
(q^{m/l} e^{-2\pi \sqrt{-1} j/l})^{\alpha} y^{m\beta} s^{m\gamma} \\
& = & \chi(X, e^{2\pi\sqrt{-1}(\tau-j)/l}, y^m)(t_1^m, t_2^m).
\end{eqnarray*}
Hence we have
\begin{eqnarray*}
&& \frac{1}{l} \sum_{j = 0}^{l-1} L((X^{lm})^{(m \sigma_l)},
    E(X^{lm}; q, y)^{(m \sigma_l)})(h_{l, j, m}, s) \\
& = &
\frac{1}{l} \sum_{j = 0}^{l-1}
\sum_{\alpha, \beta, \gamma} c(\alpha, \beta, \gamma)
(q^{m/l} e^{-2\pi \sqrt{-1} j/l})^{\alpha} y^{m\beta} s^{m \gamma} \\
& = & \sum_{\alpha, \beta, \gamma} c(l\alpha, \beta, \gamma)
q^{m\alpha} y^{m \beta} s^{m\gamma}.
\end{eqnarray*}
\end{proof}

\noindent {\em Proof of Theorem \ref{thm:Symmetric}}.
\begin{eqnarray*}
&& \sum_{n \geq 0} \chi(X^n, S_n; q, y)(s) p^n \\
& = & \sum_{n \geq 0} \sum_{\sum lN_l = n} \prod_{l \geq 1}
\frac{p^{lN_l}}{N_l! l^{N_l}}
\sum_{(N_l \sigma_l) h_l =h_l (N_l\sigma_l)}
L((X^{lN_l})^{(N_l \sigma_l)}, E(X^{lN_l}; q, y)^{(N_l \sigma_l)})(h_l, s) \\
& = & \prod_{l \geq 1} \sum_{N_l \geq 0}
\frac{1}{N_l!} \left( \frac{p^l}{l}\right)^{N_l}
 \sum_{(N_l \sigma_l) h_l =h_l (N_l \sigma_l)}
L((X^{lN_l})^{(N_l \sigma_l)}, E(X^{lN_l}; q, y)^{(N_l \sigma_l)})(h_l, s) \\
& = & \prod_{l \geq 1} C^{(l)}(s),
\end{eqnarray*}
where
\begin{eqnarray*}
&& C^{(l)}(s) \\
& = &
\sum_{N_l \geq 0}
\frac{1}{N_l!} \left( \frac{p^l}{l}\right)^{N_l}
 \sum_{(N_l \sigma_l) h_l =h_l (N_l \sigma_l)}
L((X^{lN_l})^{(N_l \sigma_l)}, E(X^{lN_l}; q, y)^{(N_l \sigma_l)})(h_l, s).
\end{eqnarray*}

Now
\begin{eqnarray*}
C^{(l)}(s) & = & \sum_{N_l \geq 0} \frac{1}{N_l!}
\left( \frac{p^l}{l}\right)^{N_l}
\sum_
{\substack{h_l = \sum N_{l, j, m} ((\sigma_l^j, 1, \dots, 1), \sigma_m)
\\
\sum_{m ,j} mN_{l,j, m} = N_l}}
\frac{l^{N_l} N_l!}{\prod_{m \geq 1} \prod_{j = 0}^{l - 1}
(lm)^{N_{l, j, m}} N_{l, j, m}!} \\
&&
L((X^{lN_l})^{(N_l \sigma_l)}, E(X^{lN_l}; q, y)^{(N_l \sigma_l)})(h_l, s)\\
& = & \sum_{N_l \geq 0} p^{lN_l}
\sum_
{\substack{h_l = \sum N_{l, j, m} ((\sigma_l^j, 1, \dots, 1), \sigma_m)
\\
\sum_{m, j} mN_{l,j, m} = N_l}}
\frac{1}{\prod_{m \geq 1} \prod_{j = 0}^{l - 1}
(lm)^{N_{l, j, m}} N_{l, j, m}!} \\
&&
\prod_{m \geq 1} \prod_{j = 0}^{l-1} \left(
L((X^{lm})^{(m \sigma_l)}, E(X^{lm}; q, y)^{(m \sigma_l)})(h_{l, j, m}, s)
\right)^{N_{l, j, m}} \\
& = & \prod_{ m \geq 1} \prod_{j=0}^{l-1} \sum_{N_{l, j, m} \geq 0}
\frac{(p^{lm})^{N_{l, j, m}}}{(l m)^{N_{l, j, m}} N_{l, j, m}!} \\
&&
\left(L((X^{lm})^{(m\sigma_l)}, E(X^{lm}; q, y)^{(m \sigma_l)})(h_{l, j, m}, s)
\right)^{N_{l, j, m}} \\
& = & \prod_{m \geq 1} \prod_{j=0}^{l-1}
\exp \left( \frac{p^{lm}}{l m}
L((X^{lm})^{(m \sigma_l)}, E(X^{lm}; q, y)^{(m \sigma_l)})(h_{l, j, m}, s)
\right) \\
& = & \exp \left( \sum_{m \geq 1} \sum_{j=0}^{l-1}\frac{p^{lm}}{l m}
L((X^{lm})^{(m \sigma_l)}, E(X^{lm}; q, y)^{(m \sigma_l)})(h_{l, j, m}, s)
\right) \\
& = & \exp \left( \sum_{m \geq 1} \frac{p^{lm}}{m}
\frac{1}{l} \sum_{j=0}^{l-1}
L((X^{lm})^{(m \sigma_l)}, E(X^{lm}; q, y)^{(m \sigma_l)})(h_{l, j, m}, s)
\right).
\end{eqnarray*}
The proof is completed by using (\ref{eqn:Lefschetz3}).

Here we have used an index calculation to prove Theorem \ref{thm:Symmetric}.
One can also prove it using the graded symmetric products generalizing the approach in \cite{Zho991, Bor-Lib}.
Such a proof is closer to a second quantized theory in flavor.


\section{Partition Functions on Framed Moduli Spaces of Instantons}
\label{sec:Framed}

In this section we recall some properties of the framed moduli spaces of instantons.
Our reference for this section is \cite{Nak-Yos}.
In particular,
we compute various equivariant genera on these spaces.
The results are expressed as sums over tuples of partitions.
Combined with the combinatorial results in \cite{Iqb-Kas1, Iqb-Kas2, Egu-Kan1, ZhoCounting, Egu-Kan2},
these will enable us to relate topological string theory with gauge theory.

\subsection{The framed moduli spaces}
Let $M(N, k)$ denote the framed moduli space of torsion free
sheaves on $\bP^2$ with rank $N$ and $c_2=k$. The framing means a
trivialization of the sheaf restricted to the line at infinity.
In particular when $N=1$ we get the Hilbert scheme $(\bC^2){[k]}$.
See \cite{Nak-Yos} for details.

As proved in \cite{Nak-Yos}, $M(N,k)$ is a nonsingular variety of dimension $2Nk$.
The maximal torus $T$ of $GL_N(\bC)$ together with the torus action on
$\bP^2$ induces an action on $M(N, k)$.
As shown in \cite{Nak-Yos}, the
fixed points are isolated and parameterized by $N$-tuples of
partitions $\vec{\mu}=(\mu^1, \cdots, \mu^N)$ such that $\sum_i |\mu^i|=k$.
The weight decomposition of the tangent bundle of $TM(N, k)$ at a
fixed point $\vec{\mu}$ is given by
\begin{eqnarray} \label{eqn:Weights}
&& \sum_{\alpha, \gamma=1}^N e_{\gamma}e_{\alpha}^{-1}
(\sum_{(i, j)\in \mu^\alpha} t_1^{-((\mu^{\gamma})^t_j - i)} t_2^{\mu_i^\alpha-j+1}
+\sum_{(i, j)\in \mu^\gamma} t_1^{(\mu^{\alpha})^t_j-i+1}t_2^{-(\mu_i^\gamma-j)}),
\end{eqnarray}
where $t_1, t_2\in \bC^*\times \bC^*,$ and $e_\alpha\in T$.

\subsection{A cohomological property of the framed moduli space}

The space $M(N, k)$ has the following remarkable property.
Let $E$ be an equivariant coherent sheaf on it.
Even though $H^i(M(N, k), E)$ might be infinite-dimensional,
the weight spaces of the induced torus action on it are finite-dimensional,
hence it makes sense to define the equivariant index.
In other words,
if
$$H^i(M(N, k), E) = \sum V_\nu$$
is the weight decomposition of $H^i(M(N, k), E)$,
then
$$\dim V_\nu < \infty$$
for all weight $\nu$.
Hence one can define
$$\ch H^i(M(N, k), E) = \sum (\dim V_\nu) e^\nu.$$
Furthermore, one can compute the equivariant index by localization
(cf. \cite{Nak-Yos}):

\begin{lemma} Let $E$ be an equivariant coherent sheaf on $M(N, k)$. Then
$$\chi(M(N, k), E)=\sum_{i=0}^{2Nk} (-1)^i \ch\, H^i(M(N, k), E)
=\sum_{\vec{\mu}}\ch
\left(\frac{i_{\vec{\mu}}^* E}{\wedge_{-1} T^*_{\vec{\mu}} M(N, k)}\right).$$
\end{lemma}

This Lemma together with (\ref{eqn:Weights}) gives us the following formula
for the equivariant elliptic genera of the framed moduli spaces:
\begin{equation} \label{eqn:Elliptic}
\begin{split}
& \sum_{k \geq 0} Q^k \chi(M(N, k), y, p)(t_1, t_2) \\
= & \sum_{\mu^{1, \dots, N}} (y^{-N}Q)^{\sum_{i=1}^N |\mu^i|}
\prod_{n \geq 1} \prod_{\alpha, \gamma} \\
& \prod_{(i, j) \in \mu^\alpha}
\frac{(1- yp^{n-1}e_{\alpha}e_{\gamma}^{-1}t_1^{(\mu^{\gamma})^t_j - i} t_2^{-(\mu_i^\alpha-j+1)})
(1- y^{-1}p^ne_{\alpha}^{-1}e_{\gamma}t_1^{-((\mu^{\gamma})^t_j - i)} t_2^{\mu_i^\alpha-j+1})}
{(1- p^{n-1} e_{\alpha}e_{\gamma}^{-1}t_1^{(\mu^{\gamma})^t_j - i} t_2^{-(\mu_i^\alpha-j+1)})
(1- p^ne_{\alpha}^{-1}e_{\gamma}t_1^{-((\mu^{\gamma})^t_j - i)} t_2^{\mu_i^\alpha-j+1})}\\
& \cdot \prod_{(i, j)\in \mu^\gamma}
\frac{(1 - y p^{n-1}e_{\alpha}e_{\gamma}^{-1}t_1^{-((\mu^{\alpha})^t_j-i+1)}t_2^{\mu_i^\gamma-j})
(1 - y^{-1}p^n e_{\alpha}^{-1}e_{\gamma}t_1^{(\mu^{\alpha})^t_j-i+1}t_2^{-(\mu_i^\gamma-j)})}
{(1 - p^{n-1}e_{\alpha}e_{\gamma}^{-1}t_1^{-((\mu^{\alpha})^t_j-i+1)}t_2^{\mu_i^\gamma-j})
(1 - p^n e_{\alpha}^{-1}e_{\gamma}t_1^{(\mu^{\alpha})^t_j-i+1}t_2^{-(\mu_i^\gamma-j)})}.
\end{split} \end{equation}
By taking $p =0$,
one can also get:
\begin{equation} \label{eqn:Chiy} \begin{split}
& \sum_{k \geq 0} Q^k \chi_y(M(N, k))(t_1, t_2) \\
= & \sum_{\mu^{1, \dots, N}} Q^{\sum_{i=1}^N |\mu^i|}
\prod_{\alpha, \gamma} \prod_{(i, j) \in \mu^\alpha}
\frac{(1- ye_{\alpha}e_{\gamma}^{-1}t_1^{(\mu^{\gamma})^t_j - i} t_2^{-(\mu_i^\alpha-j+1)})}
{(1- e_{\alpha}e_{\gamma}^{-1}t_1^{(\mu^{\gamma})^t_j - i} t_2^{-(\mu_i^\alpha-j+1)})}\\
& \cdot \prod_{(i, j)\in \mu^\gamma}
\frac{(1 - y e_{\alpha}e_{\gamma}^{-1}t_1^{-(\mu^{\alpha})^t_j-i+1)}t_2^{\mu_i^\gamma-j})}
{(1 - e_{\alpha}e_{\gamma}^{-1}t_1^{-(\mu^{\alpha})^t_j-i+1)}t_2^{\mu_i^\gamma-j})}.
\end{split} \end{equation}
If one further takes $y=0$,
\begin{equation} \label{eqn:Chi0} \begin{split}
& \sum_{k \geq 0} Q^k \chi_0(M(N, k))(t_1, t_2) \\
= & \sum_{\mu^{1, \dots, N}} Q^{\sum_{i=1}^N |\mu^i|}
\prod_{\alpha, \gamma} \prod_{(i, j) \in \mu^\alpha}
\frac{1}{(1- e_{\alpha} e_{\gamma}^{-1}t_1^{(\mu^{\gamma})^t_j - i} t_2^{-(\mu_i^\alpha-j+1)})}\\
& \cdot \prod_{(i, j)\in \mu^\gamma}
\frac{1}{(1 - e_{\alpha} e_{\gamma}^{-1}t_1^{-(\mu^{\alpha})^t_j-i+1)}t_2^{\mu_i^\gamma-j})}.
\end{split} \end{equation}

\subsection{A natural bundle on the framed moduli space}

Recall $M(N, k)$ can be identified with the space of equivalent classes of tuples of linear maps
$$(B_1: V \to V; B_2: V \to V; i: W \to V; j: V \to W)$$
satisfying
$$[B_1, B_2]+ij = 0$$
and a stability condition.
Hence one gets a vector bundle $\bV$ over $M(N, k)$ whose fibers
are given by $V$.
This bundle is an equivariant bundle,
and its weight decomposition at a fixed point $\vec{\mu}$ is given by \cite{Nak-Yos}:
\begin{eqnarray*}
&& \bV = \bigoplus_{\alpha} V_{\alpha} e_{\alpha}, \\
&& V_{\alpha} = \sum_{(i, j) \in \mu^{\alpha}} t_1^{-i+1}t_2^{-j+1}.
\end{eqnarray*}
Therefore,
the weight of ${\mbox{det}} \bV^*$ at the fixed point $\vec{\mu}$ is
$$\prod_{\alpha} \left( e_{\alpha}^{-|\mu^{\alpha}|}
\prod_{(i, j) \in \mu^{\alpha}} t_1^{i-1}t_2^{j-1}
\right).$$

Now we take $E^m_{N, k}=K_{N, k}^{\frac{1}{2}} \otimes (\det \bV^*)^m$,
where $K$ denotes the canonical line bundle of $M(N, k)$.
The equivariant index $\chi(M(N, k), E^m_{N, k})$ can be identified with
the equivariant indices of Dirac operators,
twisted by $(\det \bV^*)^m$.
An application of the above fixed point formula gives us

\begin{lemma} \label{lm:Hilbert}
We have
\begin{equation}
\begin{split}
& \sum_{k \geq 0} Q^k \chi(M(N, k), E_{N, k}^m)(e_1, \dots, e_N, t_1, t_2) \\
= & \sum_{\mu^{1, \dots, N}} Q^{\sum_{i=1}^N |\mu^i|}
\prod_{\alpha =1}^N \left( e_{\alpha}^{-|\mu^{\alpha}|}
\prod_{(i, j) \in \mu^{\alpha}} t_1^{i-1}t_2^{j-1}
\right)^m \\
& \cdot \prod_{\alpha, \gamma=1}^N  \prod_{(i, j) \in \mu^\alpha}
\frac{1}
{(e_{\alpha}^{-1}e_{\gamma}t_1^{-((\mu^{\gamma})^t_j - i)} t_2^{\mu_i^\alpha-j+1})^{\frac{1}{2}}
- (e_{\alpha}^{-1}e_{\gamma}t_1^{-(\mu^{\gamma})^t_j - i)} t_2^{\mu_i^\alpha-j+1})^{-\frac{1}{2}}} \\
& \cdot \prod_{(i, j)\in \mu^\gamma}
\frac{1}
{(e_{\alpha}^{-1}e_{\gamma} t_1^{(\mu^{\alpha})^t_j-i+1}t_2^{-(\mu_i^\gamma-j)})^{\frac{1}{2}}
 - (e_{\alpha}^{-1}e_{\gamma} t_1^{(\mu^{\alpha})^t_j-i+1}t_2^{-(\mu_i^\gamma-j)})^{-\frac{1}{2}}}.
\end{split} \end{equation}
\end{lemma}

\begin{lemma} \label{lm:Hilbert2}
When
$$t_1=e^{-\beta h},\ t_2 = e^{\beta h},\ e_\alpha = e^{-\beta a_\alpha}$$
we have the following identity:
\begin{equation*}
\begin{split}
& \sum_{k \geq 0} Q^k \chi(M(N, k), E_{N, k}^m)(e_1, \dots, e_N, t_1, t_2) \\
=& \sum_{\mu^{1, \dots, N}} (Q/4)^{\sum_{i=1}^N |\mu^i|}
\prod_{\alpha =1}^N e^{m\beta(|\mu^{\alpha}|a_{\alpha} + h \kappa_{\mu^{\alpha}}/2)} \\
& \cdot \prod_{\alpha, \gamma=1}^N
\prod_{i, j=1}^\infty \frac{\sinh{\frac{\beta}{2}}(a_{\alpha, \gamma}+h(\mu_i^\alpha
-\mu^\gamma_j+j-i))}{\sinh{\frac{\beta}{2}}(a_{\alpha, \gamma}+h(j-i))},
\end{split} \end{equation*}
where $a_{i, j} = a_i - a_j$,
\end{lemma}

\begin{proof}
By Lemma \ref{lm:Hilbert} we have
\begin{equation*}
\begin{split}
& \sum_{k \geq 0} Q^k \chi(M(N, k), E_{N, k}^m)(e_1, \dots, e_N, t_1, t_2) \\
= & \sum_{\mu^{1, \dots, N}} Q^{\sum_{i=1}^N |\mu^i|}
\prod_{\alpha =1}^N \left( e^{\beta|\mu^{\alpha}|a_{\alpha}}
\prod_{(i, j) \in \mu^{\alpha}} e^{\beta h(j-i)} \right)^m \\
& \cdot \prod_{\alpha, \gamma=1}^N  \prod_{(i, j) \in \mu^\alpha}
\frac{1}
{2\sinh \frac{\beta}{2} (a_{\alpha}-a_{\gamma} + h(\mu_i^\alpha + (\mu^{\gamma})^t_j -i-j+1))} \\
& \cdot \prod_{(i, j)\in \mu^\gamma}
\frac{1}
{2\sinh \frac{\beta}{2}(a_{\alpha} - a_{\gamma} - h((\mu^{\alpha})^t_j+\mu_i^\gamma-i-j+1))} \\
=& \sum_{\mu^{1, \dots, N}} (Q/4)^{\sum_{i=1}^N |\mu^i|}
\prod_{\alpha =1}^N e^{m\beta(|\mu^{\alpha}|a_{\alpha} + h \kappa_{\mu^{\alpha}}/2)} \\
& \cdot \prod_{\alpha, \gamma=1}^N
\prod_{i, j=1}^\infty \frac{\sinh{\frac{\beta}{2}}(a_{\alpha, \gamma}+h(\mu_i^\alpha
-\mu^\gamma_j+j-i))}{\sinh{\frac{\beta}{2}}(a_{\alpha, \gamma}+h(j-i))},
\end{split} \end{equation*}
where in the last equality we have used Lemma \ref{lm:NY} below and the identity:
$$\kappa_{\mu} = |\mu| + \sum_{i=1}^{l(\mu)} (\mu_i^2-2i\mu_i)
= 2\sum_{(i,j) \in \mu} (j-i).$$
\end{proof}

\begin{lemma} \cite{Kon} \label{lm:NY}
We have the identity:
\begin{eqnarray*}
&& \prod_{\alpha, \gamma=1}^N \left(\prod_{(i, j)\in R_\alpha}
\frac{1}{\sinh{\frac{\beta}{2}}(a_{\alpha, \gamma}+h(\mu_i^\alpha+\mu_j^{t, \gamma}-i-j+1))} \right. \\
\textbf{}&& \;\;\;\;  \cdot\left. \prod_{(i, j)\in R_\gamma}
\frac{1}{\sinh{\frac{\beta}{2}}(a_{\alpha, \gamma}-h(\mu_i^\gamma+\mu_j^{t, \alpha} -i-j+1))}\right) \\
& = & \prod_{\alpha, \gamma=1}^N
\prod_{i, j=1}^\infty \frac{\sinh{\frac{\beta}{2}}(a_{\alpha, \gamma}+h(\mu_i^\alpha
-\mu^\gamma_j+j-i))}{\sinh{\frac{\beta}{2}}(a_{\alpha, \gamma}+h(j-i))}.
\end{eqnarray*}
\end{lemma}


\section{Relationship between Hilbert Schemes and Symmetric Products}
\label{sec:Hilbert}

In this section we study the equivariant genera of $M(1, k)$ which
are the Hilbert schemes $(\bC^2)^{[k]}$.
These are given by sums over partitions by localization calculations.
By identifying them with the corresponding orbifold equivariant genera
of symmetric products of $\bC^2$,
one then gets infinite product expressions.

\subsection{Equivariant $\chi_0$ of Hilbert schemes as infinite products}
\label{sec:Chi0}

\subsubsection{Equivariant $\chi_0$ of Hilbert schemes}
By localization on $(\bC^2)^{[n]}$ (equation (\ref{eqn:Chi0})) we have
\begin{equation} \label{eqn:HilbChi0Loc}
\begin{split}
& \sum_{n=0}^\infty Q^n \chi_0((\bC^2)^{[n]})(t_1, t_2) \\
= & \sum_{n=0}^\infty Q^n \sum_{|\mu|=n}
\frac{1}{\prod_{e \in \mu} (1 - t_1^{-l(e)}t_2^{a(e)+1})(1-t_1^{l(e)+1}t_2^{-a(e)})}.
\end{split}
\end{equation}
This is a sum over partitions.

\subsubsection{Equivariant $\chi_0$ of symmetric products }
Denote by $\chi_0((\bC^2)^{(n)})(t_1, t_2)$ the orbifold equivariant geometric genera of
the symmetric product $(\bC^2)^{(n)} = (\bC^2)^n/S_n$.
Then we have
\begin{eqnarray} \label{eqn:SymmChi0Loc}
&& \chi_0((\bC^2)^{(n)})(t_1, t_2) = \exp (\sum_{n\geq 1} \frac{Q^n}{n(1-t_1^n)(1-t_2^n)}).
\end{eqnarray}
We will give three proofs.

{\em Proof by direct calculations}.
Note
$$H^i((\bC^2)^{(n)}, \cO) =\begin{cases}
\bC[x_1, y_1, \cdots , x_n, y_n]^{S_n}= S^n(\bC[x, y]), & i = 0, \\
0, & i > 0, \end{cases}$$
where $S_n$ acts on $\bC[x_1, y_1, \cdots , x_n, y_n]$
by permuting $(x_1, y_1), \cdots, (x_n, y_n)$.
Now the weight decomposition of $\bC[x,y]$ is
$$\sum_{p_1, p_2 \geq 0} t_1^{p_1}t_2^{p_2},$$
corresponding to the basis $\{x^{p_1}y^{p_2}\}$.
Then  by standard calculations for symmetric powers we have
\begin{eqnarray} \label{eqn:Hilb1}
&& \sum_{n=0}^\infty Q^n \chi_0((\bC^2)^{(n)})
= \prod_{p_1, p_2\geq 0}\frac{1}{1-t_1^{p_1}t_2^{p_2}Q}
= \exp\sum_{n=1}^\infty\frac{Q^n}{n(1-t_1^n)(1-t_2^n)}.
\end{eqnarray}

{\em Proof by orbifold equivariant localization}.
Now consider the localization on the orbifold $(\bC^2)^{(n)}$
with respect to the natural torus action.
There is only one fixed point $n [0]$ in $(\bC^2)^{(n)}$,
whose normal bundle has the following weight decomposition:
$$n (t_1 + t_2).$$
Taking into account the effect of the $S_n$-action,
one sees the contribution of this fixed point is
$$\frac{1}{n!} \frac{1}{(1-t_1)^n(1-t_2)^n}.$$

One also has to consider the twisted sectors, associated with the conjugacy classes
of $S_n$.
They can be described as follows \cite{Zho991, Zho992}.
For every partition $\mu = (\mu_1, \dots, \mu_l)$ of $n$,
let $g_{\mu}$ the product of disjoint cycles of lengths $\mu_1, \dots, \mu_l$
respectively:
$$g_{\mu} = \prod_{j=0}^{l-1} (\sum_{i=1}^j \mu_i +1,
\sum_{i=1}^j \mu_i +2, \cdots, \sum_{i=1}^j \mu_i +\mu_{j+1}).$$
The centralizer of $g_{\mu}$ is denoted by $Z_{\mu}$.
It is not hard to see that
$$Z_{\mu} \cong \prod_{i=1}^n S_{m_i(\mu)}\ltimes \bZ_i^{m_i(\mu)},$$
(semidirect product).
Hence there is a surjective homomorphism $Z_{\mu} \to \prod_{i=1}^n S_{m_i(\mu)}$,
and the order of $Z_{\mu}$ is $z_{\mu}$.

The twisted sector corresponding to $\mu$ is given by
$$(\bC^2)^{(n)}_{\mu} = ((\bC^2)^n)^{g_{\mu}}/Z_{\mu}.$$
For example,
when $\mu = (n)$,
then $g_{\mu}$ is the $n$-cylce $(1,2,\dots, n)$,
$((\bC^2)^n)^{g_{\mu}}$ is the diagonal in $(\bC^2)^n$,
$Z_{\mu}$ is the cyclic group generated by the $n$-cycle,
so the twisted sector in this case is just a copy of $\bC^2$ with a trivial $\bZ_n$-action.
There is only one fixed point situated at the origin.
The weight decomposition of its normal bundle {\em inside the whole space}
(not the twisted sector) is
$$\sum_{i=0}^{n-1} (\omega_n^i t_1 + \omega_n^it_2).$$
Here we have used the weight decomposition with respect
to the action by $T^n \times \bZ_n$.
Taking into the account of the trivial action of $\bZ_n$,
its contribution is
$$\frac{1}{n} \frac{1}{\prod_{i=0}^{n-1} (1-\omega_n^i t_1)(1-\omega_n^it_2)}
= \frac{1}{n} \frac{1}{(1-t_1^n)(1-t_2^n)}.$$

In general, it is easy to see that $((\bC^2)^n)^{g_{\mu}}$ can be identified with
$(\bC^2)^l$,
where $Z_{\mu}$ acts through the homomorphism $Z_{\mu} \to \prod_{i=1}^n S_{m_i(\mu)}$.
There is only one fixed point in each twisted sector,
given by the origin in $(\bC^2)^l$,
its normal bundle {\em inside $(\bC^2)^n$} has weight decomposition
$$\sum_{j=1}^l \sum_{i=0}^{\mu_j-1} (\omega_{\mu_j}^it_1+\omega_{\mu_j}^it_2).$$
Taking into account of the action of $Z_{\mu}$,
one gets the contribution from the twisted sector:
$$ \frac{1}{z_{\mu}} \frac{1}{\prod_{j=1}^n \prod_{i=0}^{\mu_j-1}
(1-\omega_{\mu_j}^i t_1)(1-\omega_{\mu_j}^it_2)}
= \frac{1}{z_{\mu}}\frac{1}{\prod_{i=1}^n (1-t_1^i)^{m_i(\mu)} (1-t_2^i)^{m_i(\mu)}}.$$

Therefore,
the total contribution is
\begin{eqnarray*}
&& \sum_{n \geq 0} Q^n \sum_{|\mu|=n}
\frac{1}{z_{\mu}}\frac{1}{\prod_{i=1}^n (1-t_1^i)^{m_i(\mu)} (1-t_2^i)^{m_i(\mu)}} \\
& = &  \sum_{n \geq 0} \sum_{|\mu|=n}
\prod_{i=1}^n \frac{1}{m_i(\mu)!}
\left(\frac{Q^i}{i(1-t_1^i) (1-t_2^i)}\right)^{m_i(\mu)} \\
& = & \prod_{i \geq 1} \sum_{m_i \geq 0} \frac{1}{m_i!}
\left(\frac{Q^i}{i(1-t_1^i) (1-t_2^i)}\right)^{m_i} \\
& = & \prod_{i \geq 1} \exp \left( \frac{Q^i}{i(1-t_1^i) (1-t_2^i)}\right) \\
& = & \exp \sum_{i \geq 1}\frac{Q^i}{i(1-t_1^i) (1-t_2^i)}.
\end{eqnarray*}

{\em Proof by general results for symmetric products}.
We have
\begin{eqnarray*}
&& \chi_0(\bC^2)(t_1, t_2) = \frac{1}{(1-t_1)(1-t_2)},
\end{eqnarray*}
hence by (\ref{eqn:SymmChi0}),
\begin{eqnarray*}
&& \sum_{N \geq 0} Q^N \chi_0((\bC^2)^{(n)})(t_1, t_2)
= \exp \sum_{m \geq 0} \frac{Q^m}{m}\chi_0(\bC^2)(t_1^m, t_2^m) \\
& = & \exp \sum_{m \geq 1} \frac{Q^m}{m(1-t_1^m)(1-t_2^m)}.
\end{eqnarray*}

\subsubsection{Identification of the equivariant $\chi_0$ of Hilbert schemes
and symmetric products}
We now show that
\begin{eqnarray} \label{eqn:Chi0=}
&& \chi_0((\bC^2)^{[n]})(t_1, t_2) = \chi_0((\bC^2)^{(n)})(t_1, t_2).
\end{eqnarray}
Our first proof follow \cite{Nak-Yos}.
First notice higher cohomology groups $H^i((\bC^2)^{[n]}, \cO) (i > 0)$
vanish since $(\bC^2)^{(n)}$ is a rational singularity.
Secondly,
note
$$H^0((\bC^2)^{[n]}, \cO) = H^0((\bC^2)^{(n)}, \cO),$$
hence (\ref{eqn:Chi0=}) follows.

For the second proof, note $\pi$ is a small resolution,
we have
$$\pi_*\cO_{(\bC^2)^{[n]}}  = \cO_{(\bC^2)^{(n)}}.$$
We then apply the results in \cite{Toe}.

\subsubsection{The infinite product expression for equivariant $\chi_0$ of Hilbert schemes}
Let $t_1 = q$ and $t_2 = q^{-1}$.
Putting (\ref{eqn:Chi0=}) and (\ref{eqn:SymmChi0}) together,
\begin{eqnarray*}
&& \chi_0((\bC^2)^{[n]})(q, q^{-1}) = \chi_0((\bC^2)^{(n)})(q, q^{-1}) \\
& = & \exp \sum_{m \geq 1} \frac{Q^m}{m(1-q^m)(1-q^{-m})}
= \exp \sum_{n \geq 1} \frac{- (qQ)^n}{n(1-q^n)^2}.
\end{eqnarray*}
By the following series expansion
$$\frac{1}{(1-x)^2} = \sum_{m \geq 1} m x^{m-1}.$$
we have
\begin{eqnarray*}
&& \sum_{n \geq 1} \frac{- (qQ)^n}{n(1-q^n)^2}
=   - \sum_{n\geq 1} \frac{(Q q)^n}{n}\sum_{m \geq 1} m q^{n(m-1)} \\
& = & - \sum_{ m\geq 1} m \sum_{n \geq 1} \frac{(Q q^m)^n}{n}
=  \log \frac{1}{\prod_{m \geq 1} (1 - q^mQ)^m}.
\end{eqnarray*}
Hence
\begin{eqnarray} \label{eqn:Chi0Prod}
&& \sum_{n=0}^\infty Q^n \chi_0((\bC^2)^{[n]})(q, q^{-1})
=  \frac{1}{\prod_{m \geq 1} (1 - q^mQ)^m}.
\end{eqnarray}

\subsection{Equivariant $\chi_y$ genera of Hilbert schemes as infinite products}
\label{sec:Chiy}

\subsubsection{Equivariant $\chi_y$ genera of Hilbert schemes}
By the $N=1$ case of (\ref{eqn:Chiy}) we have

\begin{equation} \label{eqn:ChiyHilb}
\begin{split}
& \sum_{n \geq 0} Q^k \chi_y((\bC^2)^{[n]})(t_1, t_2) \\
= & \sum_{\mu} Q^{|\mu|} \prod_{(i, j) \in \mu}
\frac{(1- yt_1^{-(\mu^t_j - i)} t_2^{\mu_i-j+1})(1 - y t_1^{\mu^t_j-i+1}t_2^{-(\mu_i-j)})}
{(1- t_1^{-(\mu^t_j - i)} t_2^{\mu_i-j+1})(1 - t_1^{\mu^t_j-i+1}t_2^{-(\mu_i-j)})} \\
= &  \sum_{\mu} Q^{|\mu|} \prod_{e \in \mu}
\frac{(1- yt_1^{-l(e)} t_2^{a(e)+1})(1 - y t_1^{l(e)+1}t_2^{-a(e)})}
{(1- t_1^{-l(e)} t_2^{a(e)+1})(1 - t_1^{l(e)+1}t_2^{-a(e)})}.
\end{split} \end{equation}

\subsubsection{Orbifold equivariant $\chi_y$ genera of symmetric products}
By localization formula we have
$$\chi_y(\bC^2) = \frac{(1-yt_1)(1-yt_2)}{(1-t_1)(1-t_2)}.$$
This can be directly verified as follows.
It is easy to see that
$$H^q(\bC^2, \Lambda^pT^*\bC^2)
= \begin{cases}
\bC[z_1, z_2], & (p, q) = (0, 0), \\
\bC[z_1, z_2]dz_1 \oplus \bC[z_1, z_2]dz_2, & (p, q) = (1, 0), \\
\bC[z_1, z_2]dz_1 \wedge dz_2, & (p, q) = (2, 0), \\
0, & \mbox{otherwise}.
\end{cases}
$$
So we have
\begin{eqnarray*}
\chi_y(\bC^2)(t_1, t_2)
& = & \sum_{m_1,m_2 \geq 0} (t_1^{m_1}t_2^{m_2} - y t_1^{m_1+1}t_2^{m_2}
- y t_1^{m_1} t_2^{m_2+1} + y^2 t_1^{m_1+1}t_2^{m_2+1}) \\
& = & \chi_0(\bC^2)(1-yt_1 - y t_2 + y^2t_1t_2)
= \frac{(1-yt_1)(1-yt_2)}{(1-t_1)(1-t_2)}.
\end{eqnarray*}
Therefore by standard calculations on symmetric products \cite{Zho991, Zho992} we have
\begin{eqnarray} \label{eqn:ChiySymm}
&& \sum_{n \geq 0} Q^n \chi_y((\bC^2)^{n}, S_n)(t_1, t_2)
= \exp \sum_{n \geq 1} \frac{Q^n}{n}
\frac{\chi_{y^n}(\bC^2)(t_1^n, t_2^n)}{1-y^nQ^n}.
\end{eqnarray}
Indeed,
\begin{eqnarray*}
&& \sum_{n \geq 0} Q^n \chi_y((\bC^2)^{n}, S_n)(t_1, t_2) \\
& = & \prod_{l \geq 1} \prod_{m_1,m_2 \geq 0}
\frac{(1 - y^l Q^l t_1^{m_1+1}t_2^{m_2})(1 - y^l Q^l t_1^{m_1} t_2^{m_2+1})}
{(1- Q^l y^{l-1} t_1^{m_1}t_2^{m_2}) (1 - y^{l+1} Q^l t_1^{m_1+1}t_2^{m_2+1})} \\
& = & \exp \sum_{l \geq 1} \sum_{m_1, m_2 \geq 0}
( \log (1 - y^l Q t_1^{m_1+1}t_2^{m_2})
+ \log(1 - y^l Q t_1^{m_1} t_2^{m_2+1}) \\
&& - \log(1- Q y^{l-1} t_1^{m_1}t_2^{m_2})
- \log (1 - y^{l+1} t_1^{m_1+1}t_2^{m_2+1})) \\
& = & \exp \sum_{l \geq 1} \sum_{n \geq 1} \frac{1}{n} Q^{nl} y^{n(l-1)} t_1^{nm_1}t_2^{nm_2}(1
- y^nt_1^{n} - y^nt_2^{n} + y^{2n}t_1^{n}t_2^{n}) \\
& = & \exp \sum_{n \geq 1} \frac{Q^n(1-y^nt_1^n)(1-y^nt_2^n)}{n(1- y^nQ^n)(1-t_1^n)(1-t_2^n)} \\
& = & \exp \sum_{n \geq 1} \frac{Q^n}{n}
\frac{\chi_{y^n}(\bC^2)(t_1^n, t_2^n)}{1-y^nQ^n}.
\end{eqnarray*}
This matches with (\ref{eqn:SymmChiy}).

\subsubsection{Identification of the equivariant $\chi_y$ of Hilbert schemes
and symmetric products}
Since $\pi: (\bC^2)^{[n]} \to (\bC^2)^{(n)}$ is a semismall resolution,
we have
\begin{eqnarray} \label{eqn:Chiy=}
\chi_y((\bC^2)^{[n]})(t_1, t_2) =  \chi_y (\bC^2, S_n)(t_1, t_2).
\end{eqnarray}

\subsubsection{Equivariant $\chi_y$ genera of Hilbert schemes as infinite products}
Again let $t_1 = q$, $t_2 = q^{-1}$.
By combining (\ref{eqn:Chiy=}) with (\ref{eqn:ChiySymm})
\begin{eqnarray*}
&& \sum_{n \geq 0} Q^n \chi_y((\bC^2)^{[n]})(q, q^{-1})
= \sum_{n \geq 0} Q^n \chi_y (\bC^2, S_n)(q, q^{-1}) \\
& = &\exp \sum_{n \geq 1} \frac{Q^n(1-y^nq^n)(1-y^nq^{-n})}{n(1- y^nQ^n)(1-q^n)(1-q^{-n})}.
\end{eqnarray*}
Now
\begin{eqnarray*}
&& \sum_{n \geq 1}
\frac{Q^n(1-y^nq^n)(1-y^nq^{-n})}{n(1- y^nQ^n)(1-q^n)(1-q^{-n})} \\
& = & - \sum_{n \geq 1}
\frac{Q^n(1-y^nq^n)(q^n-y^n)}{n(1- y^nQ^n)(1-q^n)^2} \\
& = & \sum_{n \geq 1} \frac{1}{n} Q^n(1-y^nq^n)(q^n-y^n)
\sum_{k \geq 0} y^{kn}Q^{kn}
\sum_{m \geq 0} (m+1)q^{nm} \\
& = & \sum_{k \geq 0} \sum_{m \geq 0}
\sum_{n \geq 1} \frac{(m+1)}{n} [ (Q^{k+1}y^kq^{m+1})^n - (Q^{k+1}y^{k+1}q^{m+2})^n \\
&& - (Q^{k+1}y^{k+1}q^{m})^n + (Q^{k+1}y^{k+2}q^{m+1})^n] \\
& = & - \sum_{k \geq 0} \sum_{m \geq 0}
(m+1)[ \log(1-Q^{k+1}y^kq^{m+1}) - \log(1-Q^{k+1}y^{k+1}q^{m+2})\\
&&- \log(1-Q^{k+1}y^{k+1}q^{m}) + \log(1-Q^{k+1}y^{k+2}q^{m+1})].
\end{eqnarray*}
Hence
\begin{equation} \label{eqn:ChiyProd}
\begin{split}
& \sum_n Q^n \chi_y((\bC^2)^{[n]}(q, q^{-1}) \\
 = & \prod_{k, m \geq 0}
\left(\frac{(1-Q^{k+1}y^{k+1}q^{m})(1-Q^{k+1}y^{k+1}q^{m+2})}
{(1-Q^{k+1}y^kq^{m+1})(1-Q^{k+1}y^{k+2}q^{m+1})}\right)^{m+1}.
\end{split} \end{equation}

\subsection{Equivariant elliptic genera of Hilbert schemes as infinite product}
\label{sec:Ell}

\subsubsection{Equivariant elliptic genera of Hilbert schemes}

The $N=1$ case of (\ref{eqn:Elliptic}) gives us
\begin{equation} \label{eqn:HilbElliptic}
\begin{split}
& \sum_{k \geq 0} Q^k \chi((\bC^2)^{[k]}; y, p)(t_1, t_2) \\
= & \sum_{\mu} (y^{-1}Q)^{|\mu|}
\prod_{n \geq 1}  \\
& \prod_{(i, j) \in \mu}
\frac{(1- yp^{n-1}t_1^{\mu^t_j - i} t_2^{-(\mu_i-j+1)})
(1- y^{-1}p^nt_1^{-(\mu^t_j - i)} t_2^{\mu_i-j+1})}
{(1- p^{n-1}t_1^{\mu^t_j - i} t_2^{-(\mu_i-j+1)})
(1- p^nt_1^{-(\mu^t_j - i)} t_2^{\mu_i-j+1})}\\
& \cdot \prod_{(i, j)\in \mu}
\frac{(1 - y p^{n-1}t_1^{-(\mu^t_j-i+1)}t_2^{\mu_i-j})
(1 - y^{-1}p^n t_1^{\mu^t_j-i+1}t_2^{-(\mu_i-j)})}
{(1 - p^{n-1} t_1^{-(\mu^t_j-i+1)}t_2^{\mu_i-j})
(1 - p^n t_1^{\mu^t_j-i+1}t_2^{-(\mu_i-j)})}.
\end{split} \end{equation}

\subsubsection{Equivariant elliptic genera of symmetric products}
By localization formula we have
\begin{equation}
\begin{split}
&  \chi(\bC^2, y, p)(t_1, t_2) \\
= & y^{-1}
\prod_{n \geq 1}
\frac{(1- yp^{n-1}t_1)(1- y^{-1}p^nt_1^{-1})}{(1- p^{n-1}t_1)(1- p^nt_1^{-1})}
\frac{(1- yp^{n-1}t_2)(1- y^{-1}p^nt_2^{-1})}{(1- p^{n-1}t_2)(1- p^nt_2^{-1})}.
\end{split} \end{equation}
By (\ref{eqn:SymmEll}) we get
\begin{eqnarray*}
& &\sum_{N \geq 0}  Q^N\chi((\bC^2)^N, S_N; p, y)(t_1, t_2) \\
& = & \exp \sum_{N, n \geq 1} \frac{Q^{Nn}}{N}
\cdot \frac{1}{n} \sum_{i=0}^{n-1} \chi(\bC^2)((\omega_n^ip^{1/n})^{N}, y^{N})(t_1^{N}, t_2^{N}).
\end{eqnarray*}

\subsubsection{Identification of the equivariant $\chi_y$ of Hilbert schemes
and symmetric products}
Motivated by \cite{Dij-Moo-Ver-Ver}
we make the following conjecture:
\begin{eqnarray}
\label{eqn:Ell=}
&& \sum_{k=0}^\infty Q^k \chi((\bC^2)^{[k]}; p, y)(t_1, t_2)
=\sum_{k=0}^\infty Q^k \chi((\bC^2)^{k}, S_k; p, y)(t_1, t_2).
\end{eqnarray}
We will refer to it as the Equivariant DMVV Conjecture.

\subsubsection{Equivariant elliptic genera of Hilbert schemes as infinite product}

Now let $t_1 = q$ and $t_2 = q^{-1}$.
Assuming the equivariant DMVV conjecture,
one then has:
\begin{eqnarray*}
&& \log \sum_{k=0}^\infty Q^k \chi((\bC^2)^{[k]}; p, y)(q, q^{-1}) \\
& = &
\log \sum_{n \geq 0} Q^n \chi((\bC^2)^n, S_n; p, y)(q, q^{-1}) \\
& = & \sum_{l \geq 1} \sum_{m \geq 1} \frac{Q^{lm}}{m}
\frac{1}{l} \sum_{j=0}^{l-1} \chi(\bC^2, e^{2\pi\sqrt{-1}(m\tau-j)/l}, y^m)(q^m, q^{-m}) \\
& = & \sum_{l \geq 1} \sum_{m \geq 1} \frac{Q^{lm}}{m}
\frac{1}{l} \sum_{j=0}^{l-1}
y^{-m} \\
&& \cdot \prod_{n \geq 1}
\frac{(1- y^me^{(n-1)2\pi\sqrt{-1}(m\tau-j)/l}q^m)(1- y^{-m}e^{n2\pi\sqrt{-1}(m\tau-j)/l}q^{-m})}
{(1- e^{(n-1)2\pi\sqrt{-1}(m\tau-j)/l}q^m)(1- e^{n2\pi\sqrt{-1}(m\tau-j)/l}q^{-m})} \\
&& \cdot \frac{(1- y^me^{(n-1)2\pi\sqrt{-1}(m\tau-j)/l}q^{-m})(1- y^{-m}e^{n2\pi\sqrt{-1}(m\tau-j)/l}q^m)}
{(1- e^{(n-1)2\pi\sqrt{-1}(m\tau-j)/l}q^{-m})(1- e^{n2\pi\sqrt{-1}(m\tau-j)/l}q^m)}\\
& = & \sum_{l \geq 1} \sum_{m \geq 1} \frac{Q^{lm}}{m}
y^{-m} \frac{1- y^mq^m}{1- q^m} \cdot \frac{1- y^mq^{-m}}{1- q^{-m}}\\
&&\cdot \frac{1}{l} \sum_{j=0}^{l-1}
\prod_{n \geq 1}
\frac{(1- y^me^{n2\pi\sqrt{-1}(m\tau-j)/l}q^m)(1- y^{-m}e^{n2\pi\sqrt{-1}(m\tau-j)/l}q^{-m})}
{(1- e^{n2\pi\sqrt{-1}(m\tau-j)/l}q^m)(1- e^{n2\pi\sqrt{-1}(m\tau-j)/l}q^{-m})} \\
&& \cdot \frac{(1- y^me^{n2\pi\sqrt{-1}(m\tau-j)/l}q^{-m})(1- y^{-m}e^{n2\pi\sqrt{-1}(m\tau-j)/l}q^m)}
{(1- e^{n2\pi\sqrt{-1}(m\tau-j)/l}q^{-m})(1- e^{n2\pi\sqrt{-1}(m\tau-j)/l}q^m)}
\end{eqnarray*}

Let
\begin{eqnarray*}
&& \prod_{n \geq 1}
\frac{(1- yp^nq)(1- y^{-1}p^nq^{-1})}
{(1- p^nq)(1- p^nq^{-1})}
\cdot \frac{(1- yp^nq^{-1})(1- y^{-1}p^nq)}{(1- p^nq^{-1})(1- p^nq)} \\
& = & \sum_{a \geq 0, b, c \in \bZ} C(a, b, c) p^aq^by^c
\end{eqnarray*}
where $C(a, b, c)$ are some integers.
Note the left-hand side is invariant under the changes of variables $q \to q^{-1}$,
or $y \to y^{-1}$,
i.e.,
\begin{eqnarray} \label{eqn:Symmetry}
&& C(a, b, c) = C(a, -b, c) = C(a, b, -c).
\end{eqnarray}
By (\ref{eqn:Symmetry}) there are integers $\tilde{C}(a, b, c)$ such that
\begin{eqnarray*}
&& \sum_{a \geq 0, b, c \in \bZ} C(a, b, c)p^aq^by^c
= \sum_{a, j \geq 0, c \in \bZ} \tilde{C}(a, j, c) p^a
(q^{2j} + q^{2(j-1)} + \cdots + q^{-2j})y^c.
\end{eqnarray*}
We have
\begin{eqnarray*}
&& \frac{1}{l} \sum_{j=0}^{l-1}
\prod_{n \geq 1}
\frac{(1- y^me^{n2\pi\sqrt{-1}(m\tau-j)/l}q^m)(1- y^{-m}e^{n2\pi\sqrt{-1}(m\tau-j)/l}q^{-m})}
{(1- e^{n2\pi\sqrt{-1}(m\tau-j)/l}q^m)(1- e^{n2\pi\sqrt{-1}(m\tau-j)/l}q^{-m})} \\
&& \cdot \frac{(1- y^me^{n2\pi\sqrt{-1}(m\tau-j)/l}q^{-m})(1- y^{-m}e^{n2\pi\sqrt{-1}(m\tau-j)/l}q^m)}
{(1- e^{n2\pi\sqrt{-1}(m\tau-j)/l}q^{-m})(1- e^{n2\pi\sqrt{-1}(m\tau-j)/l}q^m)} \\
& = &  \frac{1}{l} \sum_{j=0}^{l-1}  \sum_{a, b, c} C(a, b, c)
e^{2\pi\sqrt{-1}a(m\tau - j)/l}q^{mb}y^{mc} \\
& = &  \sum_{a, b, c} C(la, b, c)p^{ma}q^{mb}y^{mc} \\
& = & \sum_{a, j \geq 0, c \in \bZ} \tilde{C}(la, j, c) p^{ma}
\sum_{k=-j}^j q^{2km} y^{mc}.
\end{eqnarray*}
So
\begin{eqnarray*}
&& \log \sum_{k=0}^\infty Q^k \chi((\bC^2)^{[k]}; p, y)(q, q^{-1}) \\
& = & \sum_{l \geq 1} \sum_{n \geq 1} \frac{Q^{ln}}{n}
y^{-n} \frac{1- y^nq^n}{1- q^n} \cdot \frac{1- y^nq^{-n}}{1- q^{-n}} \\
&& \cdot \sum_{a, j \geq 0, c \in \bZ} \tilde{C}(la, j, c) p^{na}
\sum_{k=-j}^j q^{2km}y^{mc} \\
& = & - \sum_{n \geq 1} \frac{Q^{ln}}{n}y^{-n} (1-y^nq^n - y^nq^{-n} + y^{2n})
\sum_{m \geq 0} (m+1) q^{mn} \\
&& \cdot
\sum_{a, j \geq 0, c} \sum_{l \geq 1} \tilde{C}(la, j, c) p^{na} \sum_{k=-j}^j q^{2kn}y^{nc} \\
& = & - \sum_{m \geq 0} (m+1) \sum_{l \geq 1} \sum_{a, j \geq 0, c \in \bZ} \tilde{C}(la, j, c) \\
&& \cdot \sum_{k=-j}^j\sum_{n \geq 1} \frac{1}{n}
[(Q^lp^aq^{m+2k+1}y^{c-1})^n - (Q^lp^aq^{m+2k+2}y^{c})^n \\
&& -(Q^lp^aq^{m+2k}y^{c})^n
+ (Q^lp^aq^{m+2k+1}y^{c+2})^n] \\
& = & - \sum_{m \geq 0} (m+1)\sum_{a, j \geq 0, c} \sum_{l \geq 1} \tilde{C}(la, j, c) \\
&&\cdot
\sum_{k=-j}^j [\log(1-Q^lp^aq^{m+2k+1}y^{c-1}) - \log(1-Q^lp^aq^{m+2k+2}y^{c}) \\
&& -\log(1-Q^lp^aq^{m+2k}y^{c}) + \log(1-Q^lp^aq^{m+2k+1}y^{c+1})].
\end{eqnarray*}
Since
\begin{eqnarray*}
&& \prod_{k=-j}^j
\frac{(1-Q^2p^aq^{m+2k+2}y^{c})(1-Q^2pq^{m+2k}y^{c})}
{(1-Q^2pq^{m+2k+1}y^{c-1})(1-Q^2pq^{m+2k+1}y^{c+1})}\\
& = &
\frac{\prod_{k=-(j+\frac{1}{2})}^{j+\frac{1}{2}} (1-Q^2pq^{m+2k+1}y^{c})
\cdot \prod_{k=-(j-\frac{1}{2})}^{j-\frac{1}{2}}(1-Q^2pq^{m+2k+1}y^{c})}
{\prod_{k=-j}^j(1-Q^2pq^{m+2k+1}y^{c-1})\prod_{k=-j}^j(1-Q^2pq^{m+2k+1}y^{c+1})},
\end{eqnarray*}
therefore,
\begin{equation} \label{eqn:EllProd}
\begin{split}
& \sum_{k=0}^\infty Q^k \chi((\bC^2)^{[k]}; p, y)(q, q^{-1}) \\
= & \prod_{m \geq 0}\prod_{a, j\geq 0, c} \prod_{l \geq 1}
\left(
\frac{(1-Q^lp^aq^{b+k+2}y^{c})(1-Q^lp^aq^{b+k}y^{c})}
{(1-Q^lp^aq^{b+k+1}y^{c-1})(1-Q^lp^aq^{b+k+1}y^{c+1})}\right)^{(m+1)\tilde{C}(la, j, c)}.
\end{split}
\end{equation}

\section{Rank $1$ Examples}
\label{sec:N=1}

In this section we will study the local Calabi-Yau geometries
corresponding to gauge theory moduli spaces $M(1, k)$,
which are the Hilbert schemes $(\bC^2)^{[k]}$.
We will show that by using the results from the preceding section,
one can extract the Gopakumar-Vafa invariants
by identifying the string partition functions with the $\chi_0$, $\chi_y$
or elliptic genera of the Hilbert schemes $(\bC^2)^{[n]}$ then of
the symmetric products $(\bC^2)^{(n)}$.

\subsection{The $4D$ case}
In this subsection we consider equivariant $\chi_0$ of $M(1, k)$.
In this case the corresponding local Calabi-Yau geometry is
the resolved conifold $\cO(-1) \oplus \cO(-1) \to \bP^1$.
The Gopakumar-Vafa invariants in this case is well-known.
We rederive it from our point of view for completeness.

In this case the web diagram for the local Calabi-Yau geometry is
$$
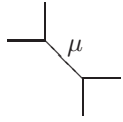
\xy
(-5, 0); (-5, 5), **@{-};
(-5, 0); (-10, 0), **@{-};
(-5, 0); (0, -5), **@{-};
(0, -5); (5, -5), **@{-};
(0, -5); (0, -10), **@{-}; (-1, -1)*+{\mu};
\endxy$$
we have following formula for the closed
string partition function \cite{Aga-Mar-Vaf, Iqb, ZhoFree}:
\begin{eqnarray} \label{eqn:Z0}
&& Z_0(Q, q) = \sum_{\mu \in \cP} (-Q)^{|\mu|} \cW_{\mu}(q)\cW_{\mu^t}(q)
= \sum_{\mu} Q^{|\mu|} \cW_{\mu}(q)\cW_{\mu}(q^{-1}).
\end{eqnarray}

We have

\begin{theorem} \label{thm:HilbChi0}
We have the identity
\begin{eqnarray} \label{eqn:HilbChi0}
&&  \sum_{\mu} Q^{|\mu|} \cW_{\mu}(t_1)\cW_{\mu}(t_2)
= \sum_{n=0}^\infty Q^n \chi_0((\bC^2)^{[n]})(t_1, t_2).
\end{eqnarray}
\end{theorem}

\begin{proof}

Recall
\begin{align*}
\cW_{\mu}(q) & = s_{\mu}(q^{\rho}), & p_{n}(q^{\rho}) & = \frac{1}{1 - q^n},
\end{align*}
where
$$q^{\rho} = (q^{-\frac{1}{2}}, q^{-\frac{3}{2}}, \dots).$$
Hence by the standard identity:
$$\sum_{\mu} s_{\mu}(x)s_{\mu}(y) = \exp \sum_{n \geq 1} \frac{p_n(x)p_n(y)}{n},$$
it follows  that
$$\sum_{\mu} Q^{|\mu|} \cW_{\mu}(t_1)\cW_{\mu}(t_2)
= \exp (\sum_{n\geq 1} \frac{Q^n}{n(1-t_1^n)(1-t_2^n)})
= \prod_{m_1, m_2 \geq 0} \frac{1}{(1 - Qt_1^{m_1}t_2^{m_2})}.$$
The proof is completed by the results in \S \ref{sec:Chi0}.
\end{proof}

\begin{remark}
By localization formula one can write the right-hand side of (\ref{eqn:HilbChi0}) as
$$\sum_{n=0}^\infty Q^n \sum_{|\mu|=n}
\frac{1}{\prod_{e \in \mu} (1 - t_1^{-l(e)}t_2^{a(e)+1})(1-t_1^{l(e)+1}t_2^{-a(e)})}.$$
Hence (\ref{eqn:HilbChi0}) is an identity that equates an infinite sum over partitions
to an infinite product:
\begin{eqnarray*}
&& \sum_{n=0}^\infty Q^n \sum_{|\mu|=n}
\frac{1}{\prod_{e \in \mu} (1 - t_1^{-l(e)}t_2^{a(e)+1})(1-t_1^{l(e)+1}t_2^{-a(e)})} \\
& = & \exp (\sum_{n\geq 1} \frac{Q^n}{n(1-t_1^n)(1-t_2^n)})
= \prod_{m_1, m_2 \geq 0} \frac{1}{(1 - Qt_1^{m_1}t_2^{m_2})}.
\end{eqnarray*}
\end{remark}

Now we take $t_1=q$ and $t_2 = q^{-1}$.
By (\ref{eqn:Z0}), (\ref{eqn:HilbChi0}) and (\ref{eqn:Chi0Prod}) we have
\begin{eqnarray*}
Z_0(Q, q) = \sum_n \chi_0((\bC^2)^{[n]})(q, q^{-1})
=  \frac{1}{\prod_{m \geq 1} (1 - q^mQ)^m}.
\end{eqnarray*}
We have a direct proof as follows.
First of all (\ref{eqn:HilbChi0Loc}) becomes
\begin{eqnarray*}
&& \sum_{n=0}^\infty Q^n \chi_0((\bC^2)^{[n]}) \\
& = & \sum_{n=0}^\infty Q^n \sum_{|\mu|=n}
\frac{1}{\prod_{e \in \mu} (1 - q^{h(e)})(1-q^{-h(e)})} \\
& = & \sum_{n = 0}^{\infty} \sum_{\mu} (-Q)^{|\mu|}
\prod_{e \in \mu} \frac{q^{h(e)/2}}{q^{h(e)/2} - q^{-h(e)/2}} \cdot
\prod_{e \in \mu} \frac{q^{-h(e)/2}}{q^{-h(e)/2} - q^{h(e)/2}} \\
& = & \sum_{n = 0}^{\infty} \sum_{\mu} (-Q)^{|\mu|} s_{\mu}(q^{\rho})s_{\mu^t}(q^{\rho}).
\end{eqnarray*}
Then the first equality follows from the identity $\cW_{\mu} = s_{\mu}(q^{\rho})$,
the second equality follows from the identity:
$$\sum_{\mu} s_{\mu}(x)s_{\mu^t}(y) = \frac{1}{\prod_{i, j} (1 + x_iy_j)},$$
where $x=(x_1, x_2, \dots)$ and $y = (y_1, y_2, \dots)$.

By comparing with (\ref{eqn:GV}) we get
$$N^j_d = \delta_{d, 1}\delta_{j, 0}.$$

\subsection{The $5D$ case}

This case has been discussed in Section 5.1.1 and Section 6.1 in \cite{Hol-Iqb-Vaf}.
We now combine their result with the results in \S \ref{sec:Chiy} to extract the Gopakumar-Vafa invariants.
The web diagram is
$$
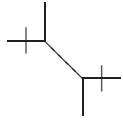
\xy
(-5, 0); (-5, 5), **@{-};
(-5, 0); (-10, 0), **@{-};
(-5, 0); (0, -5), **@{-};
(0, -5); (5, -5), **@{-};
(0, -5); (0, -10), **@{-};
(-7.5, 0)*+{|}; (2.5, -5)*+{|};
\endxy$$
The partition function given by topological vertex method is
\begin{eqnarray}
&& Z(Q, Q_m, q) = \sum_{\mu \in \cP} (-Q)^{|\mu|} \sum_{\nu \in \cP}
(-Q_m)^{|\nu|}  C_{\nu^t\mu(0)}(q)C_{\nu\mu^t(0)}(q).
\end{eqnarray}

The following result has been proved in \cite[\S 5.1.1]{Hol-Iqb-Vaf}:
\begin{eqnarray} \label{eqn:5D}
&& {\hat{Z}}=\frac{Z(Q, Q_m, q)}{Z_{(0)}(Q, q)}= \sum_{\mu } Q^{|\mu|}
\prod_{(i, j)\in \mu}
\frac{(1-Q_mq^{h(i, j)})(1-Q_mq^{-h(i,j)})}{(1-q^{h(i,j)})(1-q^{-h(i,j)})}.
\end{eqnarray}

Now we take $t_1 = q$ and $t_2 = q^{-1}$ in (\ref{eqn:ChiyHilb}) to get
\begin{eqnarray*}
&& \sum_{n \geq 0} Q^k \chi_y((\bC^2)^{[n]})(q, q^{-1})
= \sum_{\mu} Q^{|\mu|} \prod_{(i, j) \in \mu}
\frac{(1- yq^{h(i, j)})(1 - y q^{-h(i, j)})}{(1- q^{h(i, j)})(1 - q^{-h(i, j)})}.
\end{eqnarray*}
Hence we get the following result proved in \cite[\S 6.1]{Hol-Iqb-Vaf}:
After a variable change $y=Q_m$,
we have the identity
$${\hat Z} =\sum_{k=1}^\infty Q^k \chi_y(\bC^{[k]})(q, q^{-1}).$$
By (\ref{eqn:ChiyProd}) one then gets:
\begin{eqnarray} \label{eqn:ChiyGV}
&& \hat{Z} = \prod_{k, m \geq 0}
\left(\frac{(1-Q^{k+1}Q_m^{k+1}q^{m})(1-Q^{k+1}Q_m^{k+1}q^{m+2})}
{(1-Q^{k+1}Q_m^kq^{m+1})(1-Q^{k+1}Q_m^{k+2}q^{m+1})}\right)^{m+1}.
\end{eqnarray}
Comparing with (\ref{eqn:GV}) one gets:
\begin{align*}
& N_{Q^rQ_m^s}^0 = \delta_{|r-s|, 1},\\
& N_{Q^rQ_m^s}^1 = \delta_{r, s}, \\
& N_{Q^rQ_m^s}^g = 0, \;\; g > 1,
\end{align*}
where $r > 0$, $s \geq 0$.

\subsection{The $6D$ case}

This case has been discussed in Section 5.1.1 and Section 6.2 in \cite{Hol-Iqb-Vaf}.
We now combine their results with \S \ref{sec:Ell} to extract the Gopakumar-Vafa invariants in this case.
The web diagram is
$$
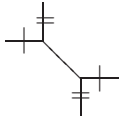
\xy
(-5, 0); (-5, 5), **@{-};
(-5, 0); (-10, 0), **@{-};
(-5, 0); (0, -5), **@{-};
(0, -5); (5, -5), **@{-};
(0, -5); (0, -10), **@{-};
(-7.5, 0)*+{|}; (2.5, -5)*+{|}; (-4.7, 2.5)*+{=}; (0, -7.5)*+{=};
\endxy$$
The partition function by topological vertex method is
\begin{eqnarray}
&& Z(Q, Q_m, Q_1, q) = \sum_{\mu, \nu, \eta \in \cP} (-Q)^{|\mu|} (-Q_1)^{|\nu|} (-Q_m)^{|\eta|}
C_{\mu\nu\eta}(q)C_{\mu^t\nu^t\eta^t}(q).
\end{eqnarray}
The following expressions for the normalized topological partition
function have been derived in \cite{Hol-Iqb-Vaf}:
\begin{eqnarray*}
&& {\hat{Z}} = \frac{Z}{Z_{(0)}} \\
& = & \sum_{\mu}Q^{|\mu|} \prod_{(i, j)\in \mu}
\frac{(1-Q_m q^{h(i, j)})(1-Q_mq^{-h(i, j)})}{(1-q^{h(i, j)})(1-q^{-h(i, j)})} \\
&& \prod_{k=1}^\infty \frac{(1-Q_\rho^kQ_m
q^{h(i, j)})(1-Q_\rho^kQ_mq^{-h(i, j)})(1-Q_\rho^kQ_m^{-1} q^{h(i, j)})
(1-Q_\rho^kQ_m^{-1}q^{-h(i, j)})}{(1-Q_\rho^kq^{h(i,j)})^2(1-Q_\rho^kq^{-h(i, j)})^2},
\end{eqnarray*}
where $Q_{\rho} = Q_1Q_m$.
The right hand side can be written in terms of the Jacobian
theta-functions.

Now let $t_1 = q$ and $t_2 = q^{-1}$ in (\ref{eqn:HilbElliptic}),
then by (\ref{eqn:HilbElliptic})
\begin{eqnarray*}
&& \sum_{k=1}^\infty Q^k \chi((\bC^2)^{[k]}; p, y)(q, q^{-1}) \\
& = & \sum_{\mu} Q^{|\mu|} \prod_{(i, j)\in \mu}
\frac{(1-yq^{h(i,j)})(1-yq^{-h(i,j)})}{(1-q^{h(i,j)})(1-q^{-h(i,j)})}\\
&& \cdot\prod_{n=1}^\infty\frac{(1-p^nyq^{h(i,j)})(1-p^ny^{-1}q^{-h(i,j)})
(1-p^nyq^{h(i,j)})(1-p^ny^{-1}q^{h(i,j)})}{(1-p^nq^{h(i,j)})^2(1-p^nq^{-h(i,j)})^2}.
\end{eqnarray*}
Thus one gets the following result in \cite[\S 6.2]{Hol-Iqb-Vaf}:
After a variable change $p=Q_\rho$ and $y = Q_m$
one has:
\begin{eqnarray*}
&& \hat{Z} = \sum_{k=1}^\infty Q^k \chi((\bC^2)^{[k]}; p, y)(q, q^{-1}).
\end{eqnarray*}
Hence if one assumes the Equivariant DMVV Conjecture,
then Gopakumar-Vafa Conjecture for this case follows by (\ref{eqn:EllProd}).
We give some examples below.

\subsection{Appendix: Examples of the 6D case}
We have
\begin{eqnarray*}
&& \sum_{a, b, c} C(a, b, c) p^aq^by^c \\
& = & 1 + p(2(q+q^{-1}) +y(q+q^{-1})+y^{-1}(q+q^{-1}))\\
&& + p^2[6-4(y+y^{-1}) + (y^2 + y^{-2})
+ (q+ q^{-1}) + 4(q^2 + q^{-2}) \\
&& - (yq+yq^{-1}+ y^{-1}q+y^{-1}q^{-1})
- 2(yq^2+yq^{-2} + y^{-1}q^2 + y^{-1}q^{-2})] + \cdots,
\end{eqnarray*}
in particular,
\begin{eqnarray} \label{eqn:C0bc}
&& C(0, b,c) =  \begin{cases}
1, & b=c=0, \\
0, & \mbox{otherwise},
\end{cases}
\end{eqnarray}
and
\begin{eqnarray} \label{eqn:C1bc}
&& C(1, b,c) =  \begin{cases}
-1, & b=\pm 1, c= \pm 1, \\
2, & b = \pm 1, c= 0, \\
0, & \mbox{otherwise}.
\end{cases}
\end{eqnarray}

\subsubsection{The $a=0$ case}

We have by (\ref{eqn:C0bc}) only one case with $C(0, b, c) \neq 0$ which gives us:
\begin{eqnarray*}
&& \prod_{k \geq 0}\prod_{b, c} \prod_{l \geq 1}
\left(
\frac{(1-Q^lq^{b+k+2}y^{c})(1-Q^lq^{b+k}y^{c})}
{(1-Q^lq^{b+k+1}y^{c-1})(1-Q^lq^{b+k+1}y^{c+1})}\right)^{(k+1)C(0, b, c)}\\
& = & \prod_{m \geq 0}\prod_{l \geq 1}
\left(
\frac{(1-Q^lq^{m+2})(1-Q^lq^{m})}
{(1-Q^lq^{m+1}y^{-1})(1-Q^lq^{m+1}y)}\right)^{(m+1)}.
\end{eqnarray*}
This matches with (\ref{eqn:ChiyGV}).

\subsubsection{The $a= 1$, $l=1$ case}
We have by (\ref{eqn:C1bc}) six cases when $C(1, b, c)$ is nonzero.
We group them into three pairs.
The pair $(b,c) = (\pm 1, 0)$ gives us:
\begin{eqnarray*}
&& \prod_{m \geq 0}
\left(
\frac{(1-Qpq^{m+3})(1-Qpq^{m+1})}
{(1-Qpq^{m+2}y^{-1})(1-Qpq^{m+2}y^{+1})}\right)^{2(m+1)} \\
&& \cdot \prod_{m \geq 0}
\left(
\frac{(1-Qpq^{m+1})(1-Qpq^{m-1})}
{(1-Qpq^{m}y^{-1})(1-Qpq^{m}y^{+1})}\right)^{2(m+1)} \\
& = & \prod_{m \geq 0}
\left(
\frac{[(1-Qpq^{m+3})(1-Qpq^{m+1})(1-Qpq^{m-1})] \cdot (1-Qpq^{m+1})}
{[(1-Qpq^{m+2}y^{-1})(1-Qpq^{m}y^{-1})][(1-Qpq^{m+2}y)(1-Qpq^{m}y)]}\right)^{2(m+1)}.
\end{eqnarray*}
The pair $(b,c) = (\pm 1, 1)$ gives us:
\begin{eqnarray*}
&& \prod_{m \geq 0}
\left(
\frac{(1-Qpq^{m+3}y)(1-Qpq^{m+1}y)}
{(1-Qpq^{m+2})(1-Qpq^{m+2}y^{2})}\right)^{-(m+1)} \\
&& \cdot \prod_{m \geq 0}
\left(
\frac{(1-Qpq^{m+1}y)(1-Qpq^{m-1}y)}
{(1-Qpq^{m})(1-Qpq^{m}y^{2})}\right)^{-(m+1)} \\
& = & \prod_{m \geq 0}
\left(
\frac{[(1-Qpq^{m+3}y)(1-Qpq^{m+1}y)(1-Qpq^{m-1}y)] \cdot (1-Qpq^{m+1}y)}
{[(1-Qpq^{m+2})(1-Qpq^{m})][(1-Qpq^{m+2}y^{2})(1-Qpq^{m}y^{2})]}\right)^{-(m+1)}.
\end{eqnarray*}
The pair $(b, c) = (\pm 1, -1)$ gives us:
\begin{eqnarray*}
&& \prod_{m \geq 0}
\left(
\frac{(1-Qpq^{m+1}y^{-1})(1-Qpq^{m-1}y^{-1})}
{(1-Qpq^{m}y^{-2})(1-Qpq^{m})}\right)^{-(m+1)} \\
&& \cdot \prod_{m \geq 0}
\left(
\frac{(1-Qpq^{m+3}y^{-1})(1-Qpq^{m+1}y^{-1})}
{(1-Qpq^{m+2}y^{-2})(1-Qpq^{m+2})}\right)^{-(m+1)} \\
& = & \prod_{m \geq 0}
\left(
\frac{[(1-Qpq^{m+3}y^{-1})(1-Qpq^{m+1}y^{-1})(1-Qpq^{m-1}y^{-1})]\cdot (1-Qpq^{m+1}y^{-1})}
{[(1-Qpq^{m+2}y^{-2})(1-Qpq^{m}y^{-2})]\cdot [(1-Qpq^{m+2})(1-Qpq^{m})]}\right)^{-(m+1)}.
\end{eqnarray*}
Note the last two pairs each contains a term
$$\prod_{m \geq 1} [(1-Qpq^{m+2})(1-Qpq^{m})]^{m+1}.$$

\subsubsection{The $a=1$ and $l=2$ case}
This involves the coefficients $C(2, b, c)$,
i.e. the coefficients of $p^2$ given as follows:
\begin{eqnarray*}
&& 6-4(y+y^{-1}) + (y^2 + y^{-2})
+ (q+ q^{-1}) + 4(q^2 + q^{-2}) \\
&& - (yq+yq^{-1}+ y^{-1}q+y^{-1}q^{-1})
- 2(yq^2+yq^{-2} + y^{-1}q^2 + y^{-1}q^{-2}) \\
&=& 2-2(y+y^{-1}) + (y^2 + y^{-2})
+ (q+ q^{-1}) + 4(q^2 + q^0+ q^{-2}) \\
&& - y(q+q^{-1})-  y^{-1}(q+q^{-1})
- 2y(q^2+q^0+q^{-2}) - 2y^{-1}(q^2 + q^{0}+q^{-2})
\end{eqnarray*}
We have rewrite the terms of the form $y^c(q^2+q^{-2})$
as
$$y^c(q^2+q^0+q^{-2}) - y^c.$$
A GV type expressions can be obtained from the term $4(q^2+q^0+q^{-2})$ as follows:
\begin{eqnarray*}
&& \prod_{m \geq 0}
\left(
\frac{(1-Q^2pq^{m+4})(1-Q^2pq^{m+2})}
{(1-Q^2pq^{m+3}y^{-1})(1-Q^2pq^{m+3}y)}\right)^{4(m+1)} \\
&& \cdot \prod_{m \geq 0}
\left(
\frac{(1-Q^2pq^{m+2})(1-Q^2pq^{m})}
{(1-Q^2pq^{m+1}y^{-1})(1-Q^2pq^{m+1}y)}\right)^{4(m+1)} \\
&& \cdot \prod_{m \geq 0}
\left(
\frac{(1-Q^2pq^{m})(1-Q^2pq^{m-2})}
{(1-Q^2pq^{m-1}y^{-1})(1-Q^2pq^{m-1}y)}\right)^{4(m+1)} \\
& = & \prod_{m \geq 0}
\left(
\frac{(1-Q^2pq^{m+1+3})(1-Q^2pq^{m+1+1})(1-Q^2pq^{m+1-1})(1-Q^2pq^{m+1-3})}
{(1-Q^2pq^{m+1+2}y^{-1})(1-Q^2pq^{m+1}y^{-1})(1-Q^2pq^{m+1-2}y^{-1})}\right)^{4(m+1)} \\
&& \cdot \prod_{m \geq 0}
\left(
\frac{(1-Q^2pq^{m+1+1})(1-Q^2pq^{m+1-1})}
{(1-Q^2pq^{m+1+2}y)(1-Q^2pq^{m+1}y)(1-Q^2pq^{m+1-2}y)}\right)^{4(m+1)}.
\end{eqnarray*}
In the same fashion one can get a GV type expression for other terms.

The case of $a=2$, $l=1$ can be dealt with similarly.
One only has to change $Q^2p$ to $Qp^2$ in the above expressions.

\section{Rank $>1$ Cases}
\label{sec:N>1}

In this section we will study local Calabi-Yau geometries
that correspond to the framed moduli spaces $M(N, k)$ for $N > 1$.
We will present details for the identifications of the relevant string partition function
as equivariant genera for suitable bundles on $M(N, k)$.
We also propose a method to obtain infinite product expressions.

\subsection{The rank $N=2$ cases}

The corresponding local Calabi-Yau geometries are given by the Hirzebruch surfaces
$\bF_m$ ($m=0, 1,2$).
Their web diagrams are given by
$$
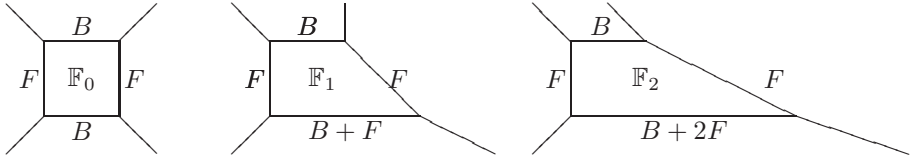
\xy
(10, 0); (5, 5), **@{-}; (15, 2)*+{B}; (45, 2)*+{B}; (84, 2)*+{B};
(10,0); (20,0), **@{-}; (25, 5); **@{-}; (8, -5)*+{F}; (78, -5)*+{F}; (38,-5)*+{F};
(10,-10); (20,-10), **@{-}; (25, -15), **@{-}; (15, -12)*+{B}; (50, -12)*+{B+F}; (95, -12)*+{B+2F};
(10,0); (10,-10), **@{-}; (5, -15), **@{-};
(20,0); (20,-10), **@{-};  (15, -5)*+{\bF_0}; (22, -5)*+{F}; (57, -5)*+{F}; (107, -5)*+{F};
(40, 0); (35, 5), **@{-}; (45, 2)*+{B};
(40, 0); (50, 0), **@{-}; (60, -10),  **@{-}; (40, -10),  **@{-};
(40, 0),  **@{-}; (47, -5)*+{\bF_1}; (38,-5)*+{F};
(40, -10); (35, -15),  **@{-};
(50, 0); (50, 5),  **@{-};
(60, -10); (70, -15),  **@{-};
(75, 5); (80, 0), **@{-}; (90, 0), **@{-}; (110, -10),  **@{-}; (80, -10),  **@{-};
(80, 0),  **@{-}; (90, -5)*+{\bF_2};
(80, -10); (75, -15),  **@{-};
(90, 0); (85, 5),  **@{-};
(110, -10); (125, -15),  **@{-};
\endxy$$
where $B, F \in H_2(\bF_m, \bZ)$ are the homological classes of the base and the fiber respectively,
and we have
\begin{align*}
F^2 & = 0, & B^2 & = -m, & BF & = 1,
\end{align*}
hence
$$(B+mF)^2 = m.$$
Their dual graphs are
$$
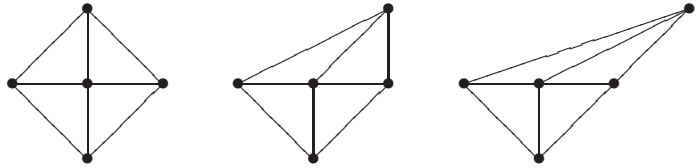
\xy
(10, -5); (20, -5),  **@{-}; (10, 5),  **@{-}; (10, -5),  **@{-};
(0, -5),  **@{-}; (10, 5), **@{-};
(0, -5); (10, -15),  **@{-}; (20, -5),  **@{-};
(10, -5); (10, -15),  **@{-};
(0, -5)*+{\bullet}; (10, -5)*+{\bullet}; (20, -5)*+{\bullet}; (10, -15)*+{\bullet};
(10, 5)*+{\bullet}; 
(40, -5); (50, -5),  **@{-}; (40, -15),  **@{-}; (40, -5),  **@{-};
(30, -5),  **@{-}; (40, -15),  **@{-};
(30, -5); (50, 5),  **@{-}; (50, -5),  **@{-};
(40, -5); (50, 5),  **@{-};
(30, -5)*+{\bullet}; (40, -5)*+{\bullet}; (50, -5)*+{\bullet}; (40, -15)*+{\bullet};
(50, 5)*+{\bullet}; 
(70, -5); (80, -5),  **@{-}; (70, -15),  **@{-}; (70, -5),  **@{-};
(60, -5),  **@{-}; (70, -15),  **@{-};
(60, -5); (90, 5),  **@{-}; (80, -5),  **@{-};
(70, -5); (90, 5),  **@{-};
(60, -5)*+{\bullet}; (70, -5)*+{\bullet}; (80, -5)*+{\bullet}; (70, -15)*+{\bullet};
(90, 5)*+{\bullet};
\endxy$$
The following results on the closed string partition functions
can be found in \cite{Aga-Mar-Vaf, Iqb, ZhoFree}:
\begin{eqnarray*}
Z_{\bF_m} & = & \sum_{\mu^{1,2,3,4}}
\cW_{\mu^1\mu^4}Q_F^{|\mu^4|} \cdot e^{-m\sqrt{-1}\kappa_{\mu^3}\lambda/2}
\cW_{\mu^4\mu^3} ((-1)^mQ_B)^{|\mu^3|} \\
&& \cdot
\cW_{\mu^3\mu^2}Q_F^{|\mu^2|} \cdot e^{m\sqrt{-1}\kappa_{\mu^1}\lambda/2}\cW_{\mu^2\mu^1}
((-1)^mQ_B)^{|\mu^1|}.
\end{eqnarray*}
Define the normalized partition function by:
\begin{eqnarray*}
{\hat{Z}}^{(m)}(Q_B, Q_F; q)
&=& \frac{Z_{\bF_m}(Q_B, Q_F; q)}{Z_{(0)}(Q_F;q)}.
\end{eqnarray*}

\begin{theorem}
After the following change of variables:
$$Q = Q_BQ_F^{-1}, e^{\beta a_1} = -1, a_1 -a_2 = 2a$$
we have the following identity:
$${\hat{Z}}^{(m)}(Q_B, Q_F; q)=\sum_{k=1}^\infty Q^k \chi(M(2,k), E^m_{2, k}(V)).$$
in other words,
the generating series of the GW invariants on the local Calabi-Yau geometry
of Hirzebruch surface $\bF_m$,
can be identified with the generating series of the equivariant indices of suitable bundles
on the framed moduli spaces.
\end{theorem}

\begin{proof}
It has been proved in \cite{Egu-Kan1, ZhoCounting} that the normalized
partition function of the topological string theory is given by
\begin{eqnarray*}
&& {\hat{Z}}^{(m)}(Q_B, Q_F; q) \\
&=& \sum_{\mu^{1, 2}}((-1)^mQ_B2^{-4}Q_F^{-1})^{|\mu^1|+|\mu^2|}Q_F^{m|\mu^2|}
q^{-\frac{m}{2}(\kappa_{\mu^1}+\kappa_{\mu^2})} \\
&& \cdot \prod_{\alpha,\gamma=1}^2
\prod_{i, j=1}^\infty
\frac{{\sinh}{\frac{\beta}{2}}(a_{\alpha,\gamma}+h(\mu_i^\alpha
-\mu^\gamma_j+j-i))}{\sinh{\frac{\beta}{2}}(a_{\alpha, \gamma}+h(j-i))},
\end{eqnarray*}
where $a_{1,2}=  -a_{2,1} = 2a$, $a_{1,1} = a_{2,2} = 0$,
$Q_F = e^{-2\beta a}$, $q = e^{-\beta h}$.
The proof is completed by Lemma \ref{lm:Hilbert2}.
\end{proof}

\subsection{The $N>2$ cases}
In this case the corresponding local Calabi-Yau geometries are
that of suitable ALE spaces of $A_{N-1}$-type fibered  over
$\bP^1$. We study the closed string partition functions of
$A_{N-1}$ fibrations.  In this case the corresponding $(p, q)$
$5$-web diagrams are of the following form:
$$
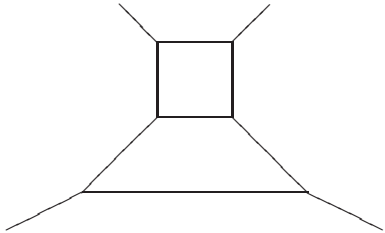
\xy
(-5, 5); (0, 0),  **@{-}; (0, -10), **@{-}; (-10, -20), **@{-}; (-20, -25), **@{-};
(15, 5); (10, 0),  **@{-}; (10, -10), **@{-}; (20, -20), **@{-}; (30, -25), **@{-};
(0, 0); (10, 0),  **@{-};
(0, -10); (10, -10),  **@{-};
(-10, -20); (20, -20),  **@{-};
\endxy$$
Their dual graph is of the form ($m=1, \dots, N+1$):
$$
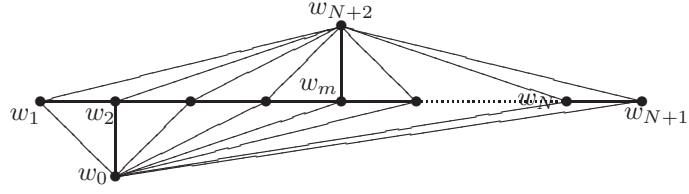
\xy
(-12, -2)*{w_1}; (-2, -2)*{w_2}; (56, 0)*{w_N}; (72, -2)*{w_{N+1}}; (27, 2)*{w_m};
(-10, 0); (0, 0), **@{-}; (10, 0), **@{-}; (20, 0), **@{-}; (30, 0), **@{-};  (40, 0), **@{-};
(60, 0), **@{.}; (70, 0), **@{-};
(-10, 0)*+{\bullet}; (0, 0)*+{\bullet}; (10, 0)*+{\bullet}; (20, 0)*+{\bullet}; (30, 0)*+{\bullet};
(40, 0)*+{\bullet}; (60, 0)*+{\bullet}; (70, 0)*+{\bullet};
(-10, 0); (0, -10), **@{-};
(0, 0); (0, -10), **@{-}; (10, 0), **@{-};
(20, 0); (0, -10), **@{-}; (30, 0), **@{-};
(40, 0); (0, -10), **@{-}; (60, 0), **@{-};
(0, -10); (70, 0), **@{-}; (0, -10)*+{\bullet};
(-10, 0); (30, 10), **@{-};
(0, 0); (30, 10), **@{-}; (10, 0), **@{-};
(20, 0); (30, 10), **@{-}; (30, 0), **@{-};
(40, 0); (30, 10), **@{-}; (60, 0), **@{-};
(30, 10); (70, 0), **@{-}; (30, 10)*{\bullet};
(-3, -10)*{w_0}; (30, 12)*{w_{N+2}};
\endxy$$

Using the topological vertex \cite{Aga-Kle-Mar-Vaf} (for a mathematical theory see \cite{LLLZ}),
one can show that \cite{Iqb-Kas2, ZhoCounting, Egu-Kan2}:
\begin{eqnarray*}
\hat{Z}^{(m)} & = & \sum_{\mu^1, \cdots, \mu^N}
((-1)^{N+m}2^{-2N}Q_B)^{|\mu^1| + \cdots + |\mu^N|} \\
&& \cdot \prod_{i=1}^{[\frac{N+m-1}{2}]} Q_{F_i}^{(N+m-2i)(|\mu^1| + \cdots + |\mu^i|)}
\prod_{i=[\frac{N+m+1}{2}]}^{N-1} Q_{F_i}^{-(N+m-2i)(|\mu^{i+1}| + \cdots + |\mu^N|)} \\
&& \cdot \prod_{i=1}^{N-1}
Q_{F_i}^{-(N-i)(|\mu^1|+\cdots +|\mu^i|)-i(|\mu^{i+1}|+\cdots +|\mu^{R_N}|)} \\
&& \cdot q^{\frac{1}{2}\sum_{i=1}^N m\kappa_{\mu^i}}\prod_{\alpha,\gamma=1}^N
\prod_{i, j=1}^\infty
\frac{{\sinh}{\frac{\beta}{2}}(a_{\alpha,\gamma}+h(\mu_i^\alpha -\mu^\gamma_j+j-i))}
{\sinh{\frac{\beta}{2}}(a_{\alpha, \gamma}+h(j-i))},
\end{eqnarray*}
where
$Q_{F_i} = e^{-\beta(a_i-a_{i+1})}$.
Here $\hat{Z}^{(m)}$ is a normalized partition function introduced in \cite{Iqb-Kas2}.

By splitting the terms with and without $m$, one gets:
\begin{eqnarray*}
\hat{Z}^{(m)} & = & \sum_{\mu^1, \cdots, \mu^N}
\left( (-1)^{N}2^{-2N}Q_B \prod_{i=1}^{[\frac{N+m-1}{2}]} Q_{F_i}^{-i}
\prod_{i=[\frac{N+m+1}{2}]}^{N-1} Q_{F_i}^{-(N-i)}\right)^{|\mu^{1}| + \cdots + |\mu^N|} \\
&& \cdot \left((-1)^{|\mu^1| + \cdots + |\mu^N|}
\prod_{i=1}^{[\frac{N+m-1}{2}]} Q_{F_i}^{|\mu^1| + \cdots + |\mu^i|}
\prod_{i=[\frac{N+m+1}{2}]}^{N-1} Q_{F_i}^{-(|\mu^{i+1}| + \cdots + |\mu^N|)}\right)^m \\
&& \cdot q^{\frac{1}{2}\sum_{i=1}^N m\kappa_{\mu^i}}\prod_{\alpha,\gamma=1}^N
\prod_{i, j=1}^\infty
\frac{{\sinh}{\frac{\beta}{2}}(a_{\alpha,\gamma}+h(\mu_i^\alpha -\mu^\gamma_j+j-i))}
{\sinh{\frac{\beta}{2}}(a_{\alpha, \gamma}+h(j-i))}.
\end{eqnarray*}

When
$$t_1=e^{-\beta h},\ t_2 = e^{\beta h},\ e_\alpha = e^{-\beta a_\alpha}$$
we have the following identity by the same method as in the proof of Lemma \ref{lm:Hilbert2}:
\begin{equation*}
\begin{split}
& \sum_{k \geq 0} Q^k \chi(M(N, k), K_{N, k}^{\frac{1}{2}} \otimes (\det \bV)^m)(e_1, \dots, e_N, t_1, t_2) \\
=& \sum_{\mu^{1, \dots, N}} (Q/2^{2N})^{\sum_{i=1}^N |\mu^i|}
\prod_{\alpha =1}^N e^{-m\beta(|\mu^{\alpha}|a_{\alpha} + h \kappa_{\mu^{\alpha}}/2)} \\
& \cdot \prod_{\alpha, \gamma=1}^N
\prod_{i, j=1}^\infty \frac{\sinh{\frac{\beta}{2}}(a_{\alpha, \gamma}+h(\mu_i^\alpha
-\mu^\gamma_j+j-i))}{\sinh{\frac{\beta}{2}}(a_{\alpha, \gamma}+h(j-i))},
\end{split} \end{equation*}
where $a_{i, j} = a_i - a_j$¡£
Thus we should take $q=e^{-\beta h}$,
\begin{eqnarray}
&& Q = (-1)^N Q_B \prod_{i=1}^{[\frac{N+m-1}{2}]} Q_{F_i}^{-i}
\prod_{i=[\frac{N+m+1}{2}]}^{N-1} Q_{F_i}^{-(N-i)},
\end{eqnarray}
and
\begin{eqnarray*}
&& \prod_{\alpha =1}^N e^{-\beta|\mu^{\alpha}|a_{\alpha}} \\
& = & (-1)^{|\mu^1| + \cdots + |\mu^N|}
\prod_{i=1}^{[\frac{N+m-1}{2}]} Q_{F_i}^{|\mu^1| + \cdots + |\mu^i|}
\prod_{i=[\frac{N+m+1}{2}]}^{N-1} Q_{F_i}^{-(|\mu^{i+1}| + \cdots + |\mu^N|)}.
\end{eqnarray*}
The right-hand side of the last equality can be rewritten as:
\begin{eqnarray*}
&& \prod_{\alpha =1}^N e^{-\beta|\mu^{\alpha}|a_{\alpha}} \cdot
\prod_{i=1}^{[\frac{N+m-1}{2}]} [(-1) e^{\beta a_k}]^{|\mu^i|} \cdot
\prod_{i=[\frac{N+m+1}{2}]}^{N} [(-1)e^{-\beta a_k}]^{|\mu^{i}|},
\end{eqnarray*}
where $k = [\frac{N+m-1}{2}]$.
Therefore,
if one takes furthermore
$$e^{\beta a_k} = 1,$$
then one has
\begin{eqnarray*}
&& \hat{Z}^{(m)} = \sum_{k \geq 0} Q^k \chi(M(N, k), K_{N, k}^{\frac{1}{2}} \otimes (\det \bV)^m)(e_1, \dots, e_N, t_1, t_2)
\end{eqnarray*}
under the above specializations.
To summarize we have:

\begin{theorem}
For $N > 2$,
we have
\begin{eqnarray}
&& \hat{Z}^{(m)} = \sum_{k \geq 0} Q^k \chi(M(N, k), K_{N, k}^{\frac{1}{2}} \otimes (\det \bV)^m)(e_1, \dots, e_N, t_1, t_2)
\end{eqnarray}
with the following specialization of variables:
$t_1=e^{-\beta h}$, $t_2 = e^{\beta h}$, $e_\alpha = e^{-\beta a_\alpha}$,
$q=e^{-\beta h}$, $e^{\beta a_{[\frac{N+m-1}{2}]}} = 1$, $Q_{F_i} = e^{-\beta(a_i-a_{i+1})}$,
and
\begin{eqnarray*}
&& Q = (-1)^N Q_B \prod_{i=1}^{[\frac{N+m-1}{2}]} Q_{F_i}^{-i}
\prod_{i=[\frac{N+m+1}{2}]}^{N-1} Q_{F_i}^{-(N-i)},
\end{eqnarray*}
\end{theorem}

\subsection{Infinite product expressions for equivariant indices}

Now we propose a method to obtain infinite product expressions for equivariant indices.

In gauge theory another partial compactification has been used.
Denote by $M_0^{\reg}(r, n)$ the framed moduli space of genuine instantons on
$S^4 = \bC^2 \cup \{\infty\}$,
and
$$M_0(r, n) = \cup_{n' = 0}^n M_0^{\reg}(r, n') \times (\bC^2)^{(n-n')}.$$
There is a projective morphism
$$\pi: M(r, n) \to M_0(r, n)$$
which is equivariant with respect the torus actions on both spaces.
In particular,
when $r=1$,
this is the Hilbert-Chow morphism:
$$\pi: (\bC^2)^{[n]} \to (\bC^2)^{(n)}.$$

We now explain why we can expect a product expression
by exploiting the collapsing morphism.
We apply the Riemann-Roch theorem in the equivariant G-theory of
algebraic stacks developed by Toen \cite{Toe}.
The rough idea is as follows.
We have converted the calculation of closed string partition functions
to the calculation of equivariant indices of $E^N_k$ on $M(N, k)$.
If we use the pushforward map $\pi_*$ in the equivariant $G$-theory,
we end up with a calculation on $M_0(N, k)$,
where localization leads to a calculation on the orbifolds $(\bC^2)^{(n)}$.
Now $G$-theory on $(\bC^2)^{(n)}$ is the same as the $S_n$-equivariant $K$-theory
of $(\bC^2)^n$.
Then by standard argument one can get a product expression
for the partition function.

\end{document}